\documentclass{amsart}

\usepackage{amsmath} 
\usepackage{amssymb}
\usepackage{mathrsfs}

\newtheorem{theorem}{Theorem}[section] 
\newtheorem{claim}[theorem]{Claim}
\newtheorem{cc}[theorem]{The Crucial Claim}

\newtheorem{observation}[theorem]{Observation}

\theoremstyle{definition}
\newtheorem{definition}[theorem]{Definition}

\newtheorem{example}[theorem]{Example}

\newtheorem{convention}[theorem]{Convention}

\newtheorem{problem}[theorem]{Problem}

\newtheorem{question}[theorem]{Question}
\newtheorem{fact}[theorem]{Fact}

\theoremstyle{remark}
\newtheorem{remark}[theorem]{Remark}
\newtheorem{notation}[theorem]{Notation}
\newtheorem{conclusion}[theorem]{Conclusion}
\newtheorem{discussion}[theorem]{Discussion}

\newcommand{\rest}{{\restriction}}
 
\newcommand{\Dom}{{\rm Dom}} 
\newcommand{\qf}{{\rm qf}} 
\newcommand{\at}{{\rm at}}

\newcommand{\Th}{{\rm Th}}

\newcommand{\LST}{{\rm LST}}
\newcommand{\wilog}{{\rm without loss of generality}}

\newcommand{\then}{{\underline{then}}}
\newcommand{\when}{{\underline{when}}}

\newcommand{\If}{{\underline{if}}}
\newcommand{\Iff}{{\underline{iff}}}
\newcommand{\mn}{{\medskip\noindent}}
\newcommand{\sn}{{\smallskip\noindent}}

\newcommand{\cL}{{\mathscr L}}

\newcommand{\bbB}{{\mathbb B}}
\newcommand{\bbI}{{\mathbb I}}

\newcommand{\gB}{{\mathfrak B}}
\newcommand{\gA}{{\mathfrak A}}

\newcommand{\cE}{{\mathscr E}}

\newcommand{\cH}{{\mathscr H}}

\newcommand{\bbL}{{\mathbb L}}
\newcommand{\bbN}{{\mathbb N}}

\newcommand{\cP}{{\mathscr P}}

\newcommand{\bbQ}{{\mathbb Q}}
\newcommand{\bbZ}{{\mathbb Z}}

\newcommand{\cU}{{\mathscr U}}

\newcommand{\cf}{{\rm cf}}

\newcount\skewfactor
\def\mathunderaccent#1#2 {\let\theaccent#1\skewfactor#2
\mathpalette\putaccentunder}
\def\putaccentunder#1#2{\oalign{$#1#2$\crcr\hidewidth
\vbox to.2ex{\hbox{$#1\skew\skewfactor\theaccent{}$}\vss}\hidewidth}}

\newenvironment{PROOF}[2][\proofname.]
   {\begin{proof}[#1]\renewcommand{\qedsymbol}{\textsquare$_{\rm #2}$}}
   {\end{proof}}

\begin{document}

\title {Nice Infinitary Logics }
\author {Saharon Shelah}
\address{Einstein Institute of Mathematics\\
Edmond J. Safra Campus, Givat Ram\\
The Hebrew University of Jerusalem\\
Jerusalem, 91904, Israel\\
 and \\
 Department of Mathematics\\
 Hill Center - Busch Campus \\ 
 Rutgers, The State University of New Jersey \\
 110 Frelinghuysen Road \\
 Piscataway, NJ 08854-8019 USA}
\email{shelah@math.huji.ac.il}
\urladdr{http://shelah.logic.at}
\thanks{The author thanks Alice Leonhardt for the beautiful typing.
The author thanks the Israel Science Foundation for partial support of
this research.  Part of this work was done while the author was
visiting Mittag-Leffler Institut, Djursholm, Sweden, in Fall 2000 and
Fall 2009.  We thank the Institut for hospitality and support.  
\newline
Publication No. 797.}

\keywords {model theory, soft model theory, characterization
  theorems, Lindstr\"om theorem, interpolation, well ordering}

\subjclass{[2010] Primary 03C95; Secondary: 03C80, 03C55}

\date{May 27, 2011}

\begin{abstract}
We deal with soft model theory of infinitary logics.  We find a logic
between $\bbL_{\infty,\aleph_0}$ and $\bbL_{\infty,\infty}$ which has
some striking properties.  First, it has interpolations (it was known
that each of those logics fail interpolation though the pair has).
Second, well ordering is not characterized in a strong way.  Third, it
can be characterized as the maximal such nice logic (in fact, is the
maximal logic
stronger than $\bbL_{\infty,\aleph_0}$ and which satisfies ``well
ordering is not characterized in a strong way").
\end{abstract}

\maketitle
\numberwithin{equation}{section}
\setcounter{section}{-1}

\centerline {Anotated Content}
\bigskip

\noindent
\S0 \quad Introduction for non-logicians, pg.2
\bigskip

\noindent
\S1 \quad Introduction and Preliminaries, pg.5
\bigskip

\noindent
\S2 \quad Introducing the logic $\bbL^1_\kappa$, pg.15
\bigskip

\noindent
\S3 \quad Serious properties of $\bbL_\kappa$, pg.22
\bigskip

\noindent
\S4 \quad $\bbL^1_\kappa$ is strong and sum/product of theories, pg.32
\newpage

\section{Introduction for non-logicians}

The first part of the introduction urged to
try to explain the aim to general mathematical audience, so
may be skipped by a knowledgeable reader; naturally we should start by
explaining what is first order logic and a general (abstract) logic
from model theoretic perspective.

We may consider classes of rings, and classes of groups but usually we
do not consider a class containing structures of both kinds.  Formally a
ring is a structure (or model) $M$ consistent of its universe, set of
elements called $|M|$ (but we may write $a \in M$) and interpretations
$+^M,\times^M$ and $0^M$ of the binary function symbols $+,\times$
and the individual consistent symbol (= zero place function symbol) $0$.
We also write ${}^n M$ for $\{(a_0,\dotsc,a_{n-1}):a_0,\dotsc,a_{n-1}$
an element of $M\}$.

Generally we have a so-called vocabulary $\tau$ consisting of relation
symbols (= predicates) and function symbol, each with a given $a$
number of places (= arity).

For a ring $M$ we many times consider the set of $n$-tuples satisfying
some equations.  Model theorist usually look at a wider class of such
sets, which start
with the family $\{\{\bar a \in {}^n M:\bar a$ satisfies an equation
$\varphi\}:n \in \bbN$ and $\varphi$ an equation$\}$ and close it
under intersection of two (with the same $n$) compliment inside the
relevant ${}^n(M)$ and projection (from ${}^{n+1}M$ to ${}^n M$).  So
a first order formula for the vocabulary $\tau,\varphi =
\varphi(x_0,\dotsc,x_{n-1})$ is a scheme giving for a $\tau$-structure
$M$ a subset $\varphi(M)$ of ${}^n M$ as above.  If $n=0,\varphi(M)
\in \{\{\langle \rangle \},\emptyset\}$ then we call 
$\varphi$ a sentence and say $M$
satisfifes it, $M \models \varphi$ iff $\varphi(M) \ne \emptyset$; let
$\bbL(\tau)$ be the set of first order sentences or formulas in the
vocabulary $\tau$; as we can add additional individual constants the
difference is minor.  We
may consider sets definable with parameters by, i.e. $\varphi(M,\bar
b) = \{\bar a:M \models \varphi[\bar a,\bar b]\}$.

There is much to be said for first order logic; e.g. for this family
of subsets of ${}^n M$ for $n \in \bbN$ has
better closure properties than the set ``solutions of finitely many
equations", i.e. varieties; however we shall not say it here.  
But in first order logic we
cannot express, e.g. ``a group $G$ is locally finite, i.e. every
finitely generated subgroup is finite".  To express this we may allow:
if $\varphi_k(x_0,\dotsc,x_{n-1})$ is a formulas for $k \in \bbN$ then so is
$\varphi(x_0,\dotsc,x_{n-1}) = \bigwedge\limits_{k} \varphi_k(x_0,\dotsc,x_k)$,
i.e. $\varphi(M) = \cap\{\varphi_k(M):k\}$.

Allowing this we get the logic $\bbL_{\aleph_1,\aleph_0}$, more
generally
\mn
\begin{enumerate}
\item[$\boxplus_0$]  the logic $\bbL_{\lambda,\kappa}$ is defined
similarly but the formulas have the form $\varphi = \varphi(\langle
x_i:i < \gamma \rangle),\gamma < \kappa$ and
we allow $\bigwedge\limits_{\alpha < \beta}
\varphi_\alpha(\bar x)$ for $\beta < \lambda$ and $(\exists
x_0,\dotsc,x_\alpha,\ldots) \varphi(\bar x,\bar y)$ where we allow a
formula to have $< \kappa$ free variables, i.e. we consider subsets of
${}^\alpha M$ for $\alpha < \kappa$.
\end{enumerate}
\mn
Another strengthening of first order logic is allowing
$\psi(\bar y) = (\exists^{\ge \aleph_1} x \varphi(x,\bar y)$,
i.e. $\psi(\bar M) = \{\bar b$: for uncountably many $a \in M$ we have
$M \models \varphi[b,\bar a]\}$.  There are many other logics.

Now first order logic has many good properties, including (recall, the
cardinality of a set $A$ is the number of elements, which may be
infinite, the cardinality of a model is the number of its element,
i.e. the cardinality of its universe $\|M\|$)
\mn
\begin{enumerate}
\item[$\boxplus_1$]  the downward \underline{LST} (L\"owenheim-Skolem-Tarski)
  property:
\sn
\begin{enumerate}
\item[$(a)$]  if a sentence $\psi \in \bbL(\tau)$ has a model, i.e. $M
  \models \psi$ then it has a countable model
\sn
\item[$(b)$]  if $\tau$ is a vocabulary, $M$ is a 
$\tau$-model, $A \subseteq |M|,|A| + |\tau|
  + \aleph_0 \le \lambda < \|M\|$ \then \, there is a $\tau$-model $N$
  of cardinality $\lambda$, a sub-model of $M$ such that $A \subseteq
  |N|$ and $\Th_{\bbL}(N) = \Th_{\bbL}(M)$ where
\sn
\item[${{}}$]  $\bullet \quad \Th_{\bbL}(M) = \{\psi \in 
\bbL(\tau_M):\psi$ a sentence $M \models \psi\}$.
\end{enumerate}
\end{enumerate}
\mn
This means that first order logic does not distinguish infinite cardinals.
\mn
\begin{enumerate}
\item[$\boxplus_2$]  compactness: if $T$ is a set of sentences in
  $\bbL(\tau)$ and every finite $T' \subseteq T$ has a model, i.e. for
  some $\tau$-model $M$ we have $\varphi \in T' \Rightarrow M \models
  \varphi$ \then \, $T$ has a model.
\end{enumerate}
\mn
The desirability of this should be obvious.
\mn
\begin{enumerate}
\item[$\boxplus_3$]  interpolation:  \If \, $\tau_0 = \tau_1
  \cap \tau_2$ are vocabularies, $\psi_1 \in \bbL(\tau_1),\psi_2 \in
  L(\tau_2)$ and $\psi_1 \vdash \psi_2$, i.e. there is no model of
  $\psi_1 \wedge \neg \psi_2$, (equivalently if $M$ is a $(\tau_1 \cup
  \tau_2)$-model and $M \models \psi_1$ then $M \models \psi$), \then
  \, there is $\varphi \in \bbL(\tau_0)$ such that $\psi_1 \vdash
  \varphi$ and $\varphi \vdash \varphi_2$.
\end{enumerate}
\mn
First order logic satisfies interpolation: this is Craig theorem.
Lindstr\"om set out to show that first order logic is the natural
choice, recalling there are many logics; for this he has first to
define a logic, essentially (see more in Definition \ref{z21})
\mn
\begin{enumerate}
\item[$\boxplus_4$]  a logic $\cL$ consists of the following
\sn
\begin{enumerate}
\item[$(a)$]  a set of sentences $\cL(\tau)$ for any vocabulary
  $\tau$, we can define formulas $\varphi(x_0,\dotsc,x_{n-1})$ by
  adding to $\tau$ individual constants
\sn
\item[$(b)$]  satisfaction relation $\models_{\cL}$, i.e. $M
  \models_{\cL} \psi$ where $M$ a model, $\psi \in \cL(\tau_M)$
\sn
\item[$(c)$]  natural properties like preservation under isomorphisms
  and monotonicity (i.e. $\tau_1 \subseteq \tau_2 \Rightarrow
\cL(\tau_1) \subseteq \cL(\tau_2))$.  
\end{enumerate}
\end{enumerate}
\mn
This seems too wide so (see more in Definition \ref{z29}).
\mn
\begin{enumerate}
\item[$\boxplus_5$]  $\cL$ is a nice logic \If \, the sets of sentences
  $\cL(\tau)$ has some natural closure properties like:
\sn
\begin{enumerate}
\item[$\bullet$]  if $\psi_1,\psi_2 \in \cL(\tau)$ then for some $\psi
  \in \cL(\tau)$ we have:
\newline
$M \models \psi$ iff $M \models \psi_1$ and $M \models \psi_2$.
\end{enumerate}
\end{enumerate}
\mn
There is a natural order on the class of logics:
\mn
\begin{enumerate}
\item[$\boxplus_6$]  $\cL_1 \le \cL_2$ \Iff \, for every vocabulary
  $\tau$ and $\psi_2 \in \cL_1(\tau)$ there is $\psi_2 \in \cL_2(\tau)$
  such that: if $M$ is a $\tau$-model then $M \models_{\cL_1} \psi_1$
\Iff \,  $M \models_{\cL_2} \psi_2$
\sn
\item[$\boxplus_7$]  $\cL_1,\cL_2$ are equivalent if $\cL_1 \le \cL_2$
  and $\cL_2 \le \cL_1$.
\end{enumerate}
\mn
Now we can phrase (really just one of the
versions\footnote{e.g. compactness just for countable theorems but
  then we have to add the occurance number is $\aleph_0$ or just
  $\aleph_0$, see Definition \ref{z26}.} of)
\mn
\begin{enumerate}
\item[$\boxtimes$]  (Lindstr\"om theorem) The logic $\cL$ is equivalent
  to $\bbL$, first order logic \when \,
\sn
\begin{enumerate}
\item[$(a)$]  $\cL$ is a nice logic
\sn
\item[$(b)$]  $\cL$ satisfies $\LST$ to $\aleph_0$:
  i.e. $\boxplus_1(a)$
\sn
\item[$(c)$]  $\cL$ satisfies compactness, see $\boxplus_2$
\end{enumerate}
\end{enumerate}
\mn
This indicates that the family of nice logics not equivalent to $\bbL$
is the union of:
\mn
\begin{enumerate}
\item[$\bullet$]  the infinitary ones, usually above
  $\bbL_{\aleph_1,\aleph_0}$
\sn
\item[$\bullet$]  the somewhat compact, usually $\aleph_0$-compact ones.
\end{enumerate}
\mn
We here deal with the first.

Lindstr\"om theorem founded ``abstract model theory" where we 
have variables over logics.  In the seventies and eighties 
this area flourished but a
reason for its almost dying out is the lack of similar theorems for
other logics, i.e. discovering (or pointing out) ``interesting" logic
which can be characterized in a reasonable way.

The aim of this work is to present such an infinitary logic and prove
that it has some desirable properties.  In particular it satifies
interpolation which holds only in ``few" cases.  This solves some more
specific old problems and we hope it will reopen the case of
``abstract model theory".

In more details, consider the logic is $\bbL^1_\kappa$ for any suitable
cardinal $\kappa$ playing the role of $\aleph_0$ in first order logic
\mn
\begin{enumerate}
\item[$\boxplus$]  $\bbL^1_\kappa$ satisfies
\sn
\begin{enumerate}
\item[$(a)$]  a downward LST; any sentence which has a model, has one
  of cardinality $< \kappa$
\sn
\item[$(b)$]  if the vocabulary has cardinality $< \kappa$ then the
  number of sentences is $\kappa$
\sn
\item[$(c)$]  a weak substitute of compactness: well ordering is not definable;
\sn
\item[$(d)$]  the $\bbL^1_\kappa$-theory of a product of two
  $\tau$-models $M_1 \times M_2$ depend just on the
  $\bbL^1_\kappa$-theories of $M_1$ and $M_2$
\sn
\item[$(e)$]  interpolation, see $\boxplus_3$.
\end{enumerate}
\end{enumerate}
\newpage

\section{Introduction and preliminaries}
\bigskip

\subsection {Aims} \
\bigskip

We feel that this is an important one among my work and will attract little
attention.  Is this an oxymoron?  We do not think so.  See below.

The investigation of model theoretic logics and soft model theory has
started with Lindstr\"om theorems and it was a central topic of model
theory in the seventies.  Major
aims were to find characterization theorems,
new important logics and non-trivial implications.  The
achievements were to a large extent summed up in the handbook
Barwise-Feferman \cite{BF} but then the subject become quite muted.  There were
some external reasons: stability theory and theoretical computer
science draw people away and there were also
some incidental personal reasons.  But probably the profound reason was
a disappointment.  The impression was that there were just too many
examples and counterexamples but not enough deep results and, 
particularly, too many logics and too
few characterization theorems (saying a logic $\bbL$ is the unique logic
such that ...).  Recall that Lindstr\"om characterize first order logic;
e.g. as the only ``reasonable" logic satisfying compactness (for
$\aleph_0$ sentences) and the downward LST theorem (for one sentence,
to $\aleph_0$), (but see \S(1D) below).
Still there was some activity later, particularly of V\"a\"an\"anen.  

Here we try to reopen the case.  A property which remains mysterious
was interpolation, see Makowsky \cite{Mw85} in the handbook.
It was known that $\bbL_{\aleph_1,\aleph_0}$ has interpolation,
(Lopez-Escobar) but not $\bbL_{\lambda,\kappa}$ when 
$(\lambda,\kappa) \ne
(\aleph_0,\aleph_0),(\aleph_1,\aleph_0)$,(Malitz).  
On $\bbL_{\kappa,\theta}$ see Dickman \cite{Dic85}.
However, the pair $(\bbL_{< \infty,\omega},
\bbL_{\infty,\infty})$ and even 
$(\bbL_{\lambda^+,\omega},\bbL_{(2^\lambda)^+,\lambda^+})$ has
interpolation, a puzzling result.  This leads naturally to a question: 
does this interpolation come from the existence of an intermediate
logic which has interpolation?  See more on the history of those
questions and on interpolation and related subjects, \cite{Mw85};

Let us recall some old questions on which we do not advance here.
Feferman raises the question
\begin{question}
\label{z6}
Is there an $\aleph_0$-compact logic strenghtening $\bbL(\exists^{\ge
\aleph_1})$ with interpolation.  
\end{question}

\noindent
Note the plethora of extensions of
$\bbL(\exists^{\ge \aleph_1})$. For my taste preferably
\begin{question}
\label{z8}
1) Is there a $\lambda$-compact logic stronger than first order
   satisfying interpolation, for any $\lambda$?

\noindent
2) Moreover, fully compact one?

Of course, part (2) becomes a question only after fully compact logics $>
\bbL$ were discovered (\cite{Sh:18}).
\end{question}

The introduction of \cite{BaFe85} mentions the (then latest advance):
some compact
logic strengthening first order logic satisfies the Beth definability
theorem, (\cite{Sh:199}) a puzzling result.  Also the pair
$(\bbL(Q^{\text{cf}}_{\aleph_0}),\bbL(aa))$ of logics satisfies interpolation.
Again a puzzling result.  Those cases give place to hope of better
results using related new logics.  Returning to infinitary logics,
old problems are (and will be our main concern):

\begin{problem}
\label{z0}
Is there a logic $\bbL$ satisfying interpolation such that $\bbL_{<
\infty,\aleph_0} \subseteq \bbL \subseteq \bbL_{< \infty,< \infty}$?
\end{problem}

\begin{problem}
\label{z1}
Assume $\kappa$ is strong limit singular of cofinality $\aleph_0$.  Is
there a logic $\bbL$ satisfying interpolation such that 
$\bbL_{\kappa^+,\aleph_0} \subseteq \bbL \subseteq \bbL_{\kappa^+,\kappa}$?
\end{problem}

Later we have asked ourselves:
\begin{problem}
\label{z2}
Is there, for arbitrarily large cardinal $\kappa$, a
logic $\bbL$ such that:
\mn
\begin{enumerate}
\item[$(a)$]  $\bbL_{\kappa,\omega} \subseteq \bbL$ and has reasonable
closure properties
\sn
\item[$(b)$]  $\bbL$ has the downward LST property in the sense
that every sentence which has a model $N$, has a model $N$ of ``small"
cardinality, moreover $M \subseteq N$
\sn
\item[$(c)$]  $\bbL$ has interpolation
\sn
\item[$(d)$]   undefinability of well ordering  (in a strong sense) which
means: if $M$ expands $(\cH(\lambda),\in)$ and $M \models \psi$ \then
\, for some $N$ we have
\begin{enumerate}
\item[$(\alpha)$]  $N \models \psi$
\sn
\item[$(\beta)$]  ord$^N$ is not well founded
\end{enumerate}

\quad a posteriori we add
\begin{enumerate}
\item[$(\gamma)$]  $|N|$ is the union of $\aleph_0$,
an internal set of bounded cardinality, i.e. for some $\langle a_n:n <
\omega\rangle,\theta,N \models ``\theta$ a cardinal such that $|a_n| \le
\theta$ and $(\forall b \in N)(\bigvee\limits_{n < \omega} N \models ``b \in
a_n")$
\end{enumerate}
\item[$(e)$]  $\bbL \subseteq \bbL_{\theta,\theta}$ for a suitable $\theta$.
\end{enumerate}
\end{problem}

\begin{problem}
\label{z4}
Is there a maximal such logic?

There is a feeling that Lindstr\"om theorem, EF-games and interpolation
are inherently connected, though I do not know of a formlization of
it, the present work gives evidence strengthening this feeling.
\bigskip

\noindent
\subsection {What is achieved} \
\bigskip

We feel that here we reasonably fulfill those old hopes mentioned
above, in the direction of non-compact logics; (recall that by
Lindstrom theorem any (nice) logic stronger than first order logic,
fail downward LST to $\aleph_0$ or fail $\aleph_0$-compactness).  

Assume for transparency $\kappa = \beth_\kappa$, we find an
interesting logic, $\bbL^1_\kappa$ such that:
\mn
\begin{enumerate}
\item[$\circledast_1$]   $(A) \quad \bbL^1_\kappa$ is a nice logic 
\sn
\item[${{}}$]  $(B) \quad$ it has a reasonable 
characterization: it is the maximal nice logic, 

\hskip25pt see Definition \ref{z21}, \ref{z29}
such that $(\alpha,<)$ can be characterized up to

\hskip25pt isomorphism  by some $\psi_\alpha \in \bbL^1_\kappa$ for
$\alpha < \kappa$, has occurance number 

\hskip25pt $\le \kappa$, see  Definition \ref{z26} and 
well ordering is not definable in a 

\hskip25pt  strong way
\sn
\item[${{}}$]  $(C) \quad$ it satisfies interpolation\footnote{for
consequences of interpolation, see \cite{Mak85}} 
(see $\boxplus_3$ above or Definition \ref{z32} below; 

\hskip25pt  answering an old question on the 
existence of such logic)
\sn
\item[${{}}$]  $(D) \quad$ it is between $\bbL^{-1}_{< \kappa} =
\cup\{\bbL_{\lambda^+,\aleph_0}:\lambda < \kappa\}$ and

\hskip25pt  $\bbL^0_{< \kappa} = 
\cup\{\bbL_{\lambda^+,\lambda^+}:\lambda < \kappa\}$, see $\boxplus_0$
above or Definition \ref{z42} below
\sn
\item[${{}}$]   $(E) \quad$ has many of the good properties of 
$\bbL^0_\kappa := \bbL^0_{< \kappa}$:
\begin{enumerate}
\item[${{}}$]  $\quad (\alpha) \quad$ downward LST, see $\boxplus_1$
above, specifically every $\psi \in \bbL^1_\kappa(\tau)$ 

\hskip35pt has a model of cardinality $< \kappa$

\hskip35pt  (see some variants in Definition \ref{a19} below)
\sn
\item[${{}}$]  $\quad (\beta) \quad$ well ordering is not an 
expressible\footnote{That is, if $\psi$ is a sentence and for every
ordinal $\alpha$ for some $M,M \models \psi$ and $(P^M,<^M)$ is a well 
ordering  of order type $\ge \alpha$ \then \, for some model $M$ of
$\psi,(P^M,<^M)$ is not a well ordering.}
\sn
\item[${{}}$]  $\quad (\gamma) \quad$ addition of theories, see \ref{a36}
\sn
\item[${{}}$]  $\quad (\delta) \quad$  product of (two) theories, see
\ref{a36}  
\end{enumerate}
\item[${{}}$]  $(F) \quad$ alternative characterization:
  $\bbL^1_\kappa$ is a minimal nice logic $\cL$ for which

\hskip35pt $\bullet \quad$ for any ordinal $\alpha < \kappa$ 
we can characterize $(\alpha,<)$ up to 

\hskip55pt isomorphism by some $\psi_\alpha \in \cL$

\hskip35pt $\bullet \quad$  we can characterize the class of $(A \cup
 \cP,A,\in)$ 

\hskip55pt where $\cP \subseteq [A]^{\le \mu}$ is an $\aleph_0$-cover
 (and $A \cap \cP = \emptyset$)

\hskip55pt  for each $\mu < \kappa$ by some sentence from $\cL$

\hskip55pt  (see Definition \ref{a24})

\hskip35pt $\bullet \quad \cL$ is $\Delta$-closed, see Definition
\ref{a31}

\hskip35pt $\bullet \quad \cL$ has occurance 
number $\le \kappa$, see Definition \ref{z26}.
\end{enumerate}
\mn
We do not have a generalization of the Feferman-Vaught theorem \cite{FeVa59}
on general operations and even not the Mostowski one, \cite{Mo52}, on
reduced products, even for the product of countably many models, see 
Theorem \ref{a43}.

Here in \S2 we choose a definition of the logic closest to the way we
arrive to it and to the proof.  For $\alpha,\theta < \kappa$ we
generalize the Ehrenfuecht-Fraisse game allowing ``rescheduling of
debts".   This does not give an
equivalence relation so we close the induced relation to an
equivalence relation and a sentence, i.e. the class of models of a
sentence is the union of some such equivalence classes; then we 
prove the basic properties.

In \S3 we deal with the deeper properties as promised in $\circledast_1$
above: non-definability of well ordering, characterization and
interpolation.

In \S4 we show how close is our logic to $\bbL^0_\kappa$ and deal with
sums and products.  We intend to continue in \cite{Sh:F1046}, in
particular concerning to \ref{z1}. 

More than once, lecturing on this some in the audience
``complain" that this definition does not sound like a definition of
logic.  So in \cite{Sh:F1046} we intend to give presentation close to the
ways logic are traditionally defined 
(we could have done it for $\bbL^1_\kappa$, too) \underline{but} our
characterization theorem shows that we shall get the same logic.

The logic $\bbL^1_\kappa$ from \S1 is quite satisfactory; many of the
good properties of $\bbL_{< \kappa,\aleph_0}$ and interpolation, and a
characterization (parallel to Lindstr\"om theorem).  \underline{But}
compared to $\bbL_{< \kappa,\aleph_0}$ we lose the upward LST.
\end{problem}

\begin{question}
\label{z9}
Let $\kappa = \beth_\kappa$, is there $\bbL$ such that $\bbL^{-1}_\kappa
\le \cL \le \bbL^0_\kappa$
satisfying interpolation and the upward LST theorem?

We intend to deal with this in \cite{Sh:F1046}.
\end{question}
\bigskip

\subsection {Why characterizations?}

Note that characterization theorems are central for several reasons:
\mn
\begin{enumerate}
\item[$\circledast_2$]  $(a) \quad$ per se, uniqueness 
results are nice, of course
\sn
\item[${{}}$]  $(b) \quad$ historically - Lindstr\"om theorem has this form
\sn
\item[${{}}$]  $(c) \quad$ they prove a logic is natural logic
\sn
\item[${{}}$]  $(d) \quad$ for a logic which lacks such theorem
we may well suspect that
\begin{enumerate}
\item[${{}}$]  $\bullet \quad$ there are many relatives of similar 
good properties, without a 

\hskip35pt special reason to prefer one or another.  
\end{enumerate}
\end{enumerate}
\mn
How good is a characterization theorem?
Of course, it all depends on the properties appearing in the
characterization being natural and preferably well established.

The situation of having reasonable logics which we can strengthen
preseving their main positive properties but neither seeing a maximal
one nor proving such extension do not exist has been prominent in the
area; e.g. the most well established ones like $\bbL(Q) = \bbL(\exists^{\ge
\aleph_1})$ and also $\bbL_{\kappa^+,\aleph_0}$.

Why here we tend to look at strong limit cardinals, in particular,
$\kappa = \beth_\kappa$?  Note that in first order formula with a fix
finite vocabulary, with predicates only for transparency 
the number of sentences of quantifier depth $q$ has order of magnitude
$\beth_q$, iterated power $q$ times.  In infinitary logic,
if we like that still there is a sentence expressing the ``quantifier
depth $\le \alpha$ theory" we need $\kappa = \beth_\kappa$.  A price is
that for $\kappa$ singular, we lose full substitution, and moreover,
full closure under conjunctions of $< \kappa$.  This is resolved if 
we demand $\kappa$
is strong limit regular, i.e. (strongly) inaccessible, fine but the
existence of such cardinals is unprovable in ZFC.
\bigskip

\subsection {Directions ignored here}

We do not deal here with some other major directions:
\mn
\begin{enumerate}
\item[$\circledast_3$]  $(a) \quad \aleph_0$-compact (nice logic)
\sn
\item[${{}}$]  $(b) \quad$ logics without negation and continuous logic
\sn
\item[${{}}$]   $(c) \quad$ almost isomorphism (and absolute logic).
\end{enumerate}
\mn
We may look at model theory essentially replacing ``isomorphic" by
``almost isomorphic", that is isomorphisms by
potential isomorphisms, i.e. isomorphism in some forcing extension.  In
\cite{Sh:12} we have suggested to reconsider a major theme in model
theory, counting the number of isomorphism types.  We call $M,N$
almost-isomorphic when $M,N$ have (the same vocabulary and) the same
$\bbL_{\infty,\aleph_0}$-theory, equivalently isomorphic in some
generic extension.  For a theory $T$ let
$\bbI_{\text{ai}}(\lambda,T)$ be
$|\{M/\equiv_{\bbL_{\infty,\aleph_0}}:M$ a model of $T$ of cardinality
$\lambda\}|$.  This behaves nicely: if $T$ has cardinality $\le
\lambda$, is first order or just $\subseteq \bbL_{\lambda^+,\aleph_0}$ then
$\dot{\bbI}_{\text{ai}}(\lambda,\psi) 
\le \lambda < \mu \Rightarrow \dot{\bbI}_{\text{ai}}(\mu,T) \le 
\dot{\bbI}_{\text{ai}}(\lambda,T)$, (on $\dot{\bbI}_{a_i}(-,T)$ for
$\aleph_0$-stable $T$, see a work of Laskowski-Shelah in
preparation).  In \cite{Sh:12} we also define ``$M$ is ai-rigid,
i.e. $a \ne b \in M \Rightarrow (M,a) \not\equiv_{\bbL_{\infty,\aleph_0}}
(N,a)"$ and have downward LST theorem for it.   Later
Nadel suggested further to consider homomorphisms, in particular for
abelian groups, see \cite[Ch.IV,\S3,pg.487]{EM02}, more G\"obel-Shelah
\cite{GbSh:880}, G\"obel-Herden-Shelah \cite{GbHeSh:948}.
Barwise characterized the relevant logic, $\bbL_{\infty,\aleph_0}$ by
absoluteness: (among logics with satisfaction being absolute under
forcing it is maximal).
\bigskip

\subsection {Can soft model theory be applied?} \
\bigskip

What about applications of soft model theory?  There are some
applications using compact logics.  
See \cite{Sh:384} on extending first order logic by second order
quantifiers restricted in some ways.  Of course, the expressive power
of the logic depend on the restriction, see the example below.  In
cofinality logic, $\bbL(\bold Q^{\text{cf}}_C)$ with $\bold C$ is a
class of regular cardinals, we are allowed to say: the formula
$\varphi(x,y)$, possibly with parameters, define a linear order with
no last element of cofinality from $\bold C$, recalling that the
cofinality of a linear oder $I$ is the minimal cardinality of an
unbounded (equivalently cofinal) subset.  If $\bold C$ is
non-trivial, this logic is
a very interesting logic (e.g. fully compact), in
particular showing what we cannot prove so full compactness is not
sufficient to characterize first order logic.  But its expressive power
is weak so we do not expect it to have applications.

In \cite{Sh:384} (where you can find something about the history of
the topic) we prove the compactness of the quantifier (that is for first
order logic extended by it) ${\cL}^{\text{\rm ceab}}=
\bbL(\bbQ^{\text{\rm ceab}})$ --- quantifying over complete 
embeddings of one atomless Boolean ring into another.  Moreover, for
this logic we prove completeness for a natural set of axioms.  Now consider
the problem ``can the automorphism groups of a
1-homogeneous\footnote{A Boolean algebra $B$ is 1-homogeneous if it is
atomless for every $a$, $b\in B \setminus \{0_{\bbB}\}$ we have 
$B \cong B \restriction b$ (equivalently for $a$,
$b\in \bbB \setminus \{0_B,1_B\}$ for some automorphism $f$ of $B$,
$f(a)=b$)}
Boolean algebra be non-simple\footnote{That is has no normal
subgroup which is neither the full group nor the one-element sub-group}"?
Much is known on this group and, in particular, that it 
is ``almost" simple --- see Rubin-Stepanek \cite{RS89}.  It
was known that there may exist such Boolean Algebras as by 
\cite[IV]{Sh:b} in some generic extension, all automorphisms of  
${\cP}(\omega)/ \text{\rm finite}$ are trivial (i.e. induced by
permutations $\pi$ of $\bbZ$ such that $\{n \in \bbZ:n \ge 0$ but 
$\pi(n) < 0\}$ is finite) and van Dowen note that the 
group of trivial automorphisms of
${\cP}(\omega)/\text{\rm finite}$ is not simple 
(as the subgoup of the automorphisms induced by permutations of
$\omega$ is a normal subgroup) and the quotient is isomorphic to
$(\bbZ,+)$.  Alternatively,
Koppelberg \cite{Kp85} has directly constructed such Boolean Algebras of
cardinality $\aleph_{1}$ assuming (the more natural assumption) CH.  So by the
completeness theorem (as the set of axioms is absolute), as the relevant facts 
are expressible in ${\bbL}(\bbQ^{\text{\rm ceab}})$, 
the existence is proved in  ZFC.  Some may want to consider a
direct proof.  It almost certainly will give more specific desirable
information.

Another helpful quantifier is on branches of trees (see
\cite{Sh:72}).  In Fuchs-Shelah \cite{FuSh:766} it is used to 
eliminate the use of diamond, i.e. to prove in ZFC the existence of
valuation domains $R$ such that there are $R$-modules which are
univerasl (i.e. the family of sub-modules is linearly ordered by
inclusion) but not standard.  Note the obvious examples (which are
called standard): $R$ itself or
appropriate quotients.  Actually the completeness theorem for this
logic gives an absoluteness result which is used.  

We believe that quantifiers with completeness and 
compactness will be useful so it is worthwhile to find such
quantifiers.  Hopefully see more in \cite{Sh:800}.
\bigskip

\centerline {$* \qquad * \qquad *$}
\bigskip

\subsection {Preliminaries} \
\bigskip

\begin{notation}
\label{z10}
1) $\tau$ denotes a vocabulary, i.e. a set of predicates so each $P
   \in \tau$ has arity$_\tau(P) < \omega$ places and each function
symbol $F \in \tau$ has arity$_\tau(F) < \omega$ places; of
 course, individual constants are zero-place function symbols and we
   may write arity$(P)$, arity$(F)$ when $\tau$ is clear from the context.

\noindent
1A) For a structure $M$ let $\tau_M$ be the vocabulary of $M$, for a
predicate $P$ from $\tau_M,P^M$ is the interpretation of $P$ so an
arity$_\tau(P)$-place relation on $|M|$, the universe of $M$;
similarly for a function symbol $F$ from $\tau$ and in particular for
an individual constant $c$ from $\tau$.

\noindent
2) $\cL$ denotes a logic, see Definition \ref{z21}.

\noindent
3) $\bar x,\bar y,\bar z$ denote sequence of variables (with no
   repetition).  Usually $\bar x = \langle x_i:i < \alpha\rangle$ so
   $\alpha = \ell g(\bar \alpha)$ but even possibly $\bar x = \langle x_s:s
\in S\rangle$ and then we let $\ell g(\bar x) = S$.

\noindent
4) We say $\tau$ is a relational vocabulary \when \, it has no
function symbol.
\end{notation}

\noindent
Recall
\begin{definition}
\label{z21}
1) A logic $\cL$ consists of
\mn
\begin{enumerate}
\item[$(a)$]  function $\tau \mapsto \cL(\tau)$ giving a set of
sentences $\varphi$ (or formulas $\varphi(\bar x)$, see \ref{z24} below) for
any vocabulary $\tau$; the function is a class function that is a definition
\sn
\item[$(b)$]  $\models_{\cL}$, satisfaction, i.e. the relation $M \models_{\cL}
\varphi$ for $M$ a model, $\varphi \in \cL(\tau_M)$
\sn
\item[$(c)$]  renaming: the function $\hat \pi$, depending on
$(\tau_1,\tau_2,\pi,\cL)$, is a one-to-one function from $\cL(\tau_1)$
onto $\cL(\tau_2)$ \when \, $\pi$ is an isomorphism from the vocabulary
$\tau_1$ onto the vocabulary $\tau_2$ (i.e. if $P \in \tau_1
\Rightarrow$ then $\pi(P) \in \tau_2$ is a predicate and
arity$_{\tau_1}(P) = \text{ arity}_{\tau_2}(\pi(P))$ and similarly for
$F \in \tau_1$)
\sn
\item[$(d)$]  if $\pi$ is an isomorphism from the vocabulary
$\tau_1$ onto the vocabulary $\tau_2$ and $M_1$ is a $\tau_1$-model and
$M_2 = \pi(M_1)$ naturally defined \then \, $\varphi \in
\cL(\tau_1(M_1)) \Rightarrow [M_1 \models \varphi \Leftrightarrow
M_2 \models \hat \pi(\varphi)]$
\sn
\item[$(e)$]  (isomorphism): if $M_1,M_2$ are isomorphic $\tau$-models
and $\varphi \in \cL(\tau)$ then $M_1 \models_{\cL} \varphi
\Leftrightarrow M_2 \models_{\cL} \varphi$
\sn
\item[$(f)$]  (monotonicity): if $\tau_1 \subseteq \tau_2$ \then \,
$\cL(\tau_1) \subseteq \cL(\tau_2)$ and for any $\tau_2$-model $M_2$ and
$\varphi \in \cL(\tau_1)$ we have $M_2 \models \varphi
\Leftrightarrow (M_2 \rest \tau_1) \models \varphi$.
\end{enumerate}
\end{definition}

\begin{convention}
\label{z24}
We define a formula $\varphi = \varphi(\bar x)$ in $\cL(\tau)$ as a
sentence in $\cL(\tau \cup\{c_i:i < \ell g(\bar x)\})$ with $c_i(i <
\alpha)$ pairwise distinct individual constants not from $\tau$
 and if
$\varphi(\bar x) \in \cL(\tau)$ and $\bar a \in {}^{\ell g(\bar x)}M$
then $M \models_{\cL} \varphi[\bar a]$ means $M^+ \models_\cL
\varphi$ where $M^+$ is the expansion of $M$ by $c^{M^+}_i = a_i$ for
$i < \ell g(\bar x)$, in fact, $\ell g(\bar x)$ can be any index set.
\end{convention}

\begin{definition}
\label{z26}
For a logic $\cL$, the occurance number oc$(\cL)$ of $\cL$ is the
 minimal cardinal $\kappa$ (or $\infty$) such that $\cL(\tau) =
   \cup\{\cL(\tau'):\tau' \subseteq \tau$ is of cardinality $<
   \kappa\}$ for any vocabulary $\tau$.
\end{definition}

\begin{definition}
\label{z29}
We say a logic $\cL$ is nice \when :
\mn
\begin{enumerate}
\item[$(a)$]   $(\alpha) \quad$ applying predicates:  if $P \in \tau$ 
is an $n$-place predicate then 

\hskip25pt $P(x_0,\dotsc,x_{n-1}) 
\in \cL(\tau)$ is a formula of $\cL(\tau)$
\sn
\item[${{}}$]    $(\beta) \quad$ equality: $x_0 = x_1 \in \cL(\tau)$, i.e.
for any individual constants $c_0,c_1 \in \tau$

\hskip25pt  there is $\varphi
\in \cL(\tau)$ such that $M \models \varphi$ iff $c^M_0 = c^M_1$
\sn
\item[${{}}$]   $(\gamma) \quad$ applying\footnote{We may put together
 clauses $(a)(\alpha)$ and $(a)(\gamma)$ allowing $p(\sigma_0(\bar
 x),\dotsc,\sigma_{n-1}(\bar x)),\sigma_\ell$ a term; the choice is
 immaterial.}  functions: for an $n$-place
function symbol $F \in \tau$,

\hskip25pt $x_0 = F(x_1,\dotsc,x_{n-1})$ is a
formula of $\cL(\tau)$
\sn
\item[$(b)$]   $(\alpha) \quad \cL$ is closed under conjunction, i.e. for every
$\varphi_1,\varphi_2 \in \cL(\tau)$ there is 

\hskip25pt $\varphi_3 \in \cL(\tau)$
such that for every $\tau$-model $M$ we have:

\hskip25pt  $M \models \varphi_3$ \Iff \,
$M \models \varphi_1$ and $M \models \varphi_2$
\sn
\item[${{}}$]   $(\beta) \quad \cL$ is closed under existential quantifier,
$(\exists x)$, i.e. for any $\varphi \in \cL(\tau \cup$

\hskip25pt $\{c\}),c$ an
individual constant not in $\tau$, there is 
$\psi \in \cL(\tau)$ such

\hskip25pt  that for any
$\tau$-model $M$ we have $M \models \varphi$ iff $M^+ \models \varphi$
for some 

\hskip25pt $(\tau \cup \{c\})$-expansion $M^+$ of $\mu$
\sn
\item[${{}}$]   $(\gamma) \quad \cL$ is closed under negation, i.e. 
for any $\varphi \in \cL(\tau)$ there is $\psi \in \cL(\tau)$ 

\hskip25pt such that for any $\tau$-model 
$M,M \nVdash_{\cL} \varphi \Leftrightarrow M \models_{\cL} \psi$
\sn
\item[$(c)$]  restricting a sentence $\psi$ to a predicate $P,\psi
\rest P \in \cL(\tau)$ when $\psi \in \cL(\tau)$, 
where $P \in \tau$ is a unary predicate, and $\tau$ is a
relational vocabulary, see Definition \ref{z34}, \ref{z10}(4),
(but see clause $(d)(\gamma)$ below)
\sn
\item[$(d)$]  weak substitution, that is, has substitution for very
  simple schemes, see \ref{z37} below.
\end{enumerate}
\end{definition}

\begin{remark}
\label{z31}
1) Above we prefer $(a)(\gamma)$ on using 
$R'(\sigma_0,\dotsc,\sigma_{m-1})$
where each $\sigma_\ell = \sigma_\ell(\bar x)$ is a term.

\noindent
2) Below we can define the multi-sort vesrion.
\end{remark}

\begin{definition}
\label{z32}
1) A logic $\cL$ satisfies interpolation \when \, for any sentence
   $\psi_1 \in \cL(\tau_1),\varphi_2 \in \cL(\tau_2)$ and $\tau =
   \tau_1 \cap \tau$, we have: if $\varphi_1 \vdash \varphi_2$,
   i.e. for any $(\tau_1 \cup \tau_2)$-model $M, M \models \varphi_1
   \Rightarrow M \ne \varphi_2$, \then \, for some sentence 
$\psi \in \cL(\tau)$ we have $\varphi_1 \vdash \psi$ and $\psi \vdash
   \varphi_2$.

\noindent
2) Naturally definition for multi-sort languages.
\end{definition}

\begin{definition}
\label{z34}
1) We say $\varphi \equiv \psi \rest P$ where $\varphi,\psi \in
\cL(\tau)$ or pedantically\footnote{of course $\psi \rest \tau$ is not
uniquely determined.  We may like to be more liberal in restricting a model,
in \ref{z10}(1) allow $F^M$ to be a partial function for $F$ a function
symbol from $\tau(M)$.  We may combine this with restriction (see Definition
\ref{z34}), then we still have to demand $P^M \ne \emptyset$
(except if we go further and allow empty models).  Note $M \models
``\neg F(\bar a) = b"$ when $F^M(\bar a)$ is not well defined.}
$\varphi = \psi
\rest_\tau P$ \underline{where} $P$ is a unary predicate in the
vocabulary $\tau,\tau$ is of minimal cardinality such
that $\psi \in \cL(\tau)$ \when \,: 

for any $\tau$-models $M,M \models
(\psi \rest_\tau P)$ \Iff \, $P^M$ is non-empty, closed under $F^M$ for
every function symbol from $\tau$ and $M \rest P^M \models \psi$, see below.

\noindent
2) For a $\tau$-model $M$ and unary predicate $P$ let $N = M \rest
   \tau$ be the $\tau$-model with universe $P^M$, for $n$-place
   predicate $Q \in \tau$ we have $Q^N = Q^M \cap {}^n|N|$ and for
   $n$-place function symbol $F \in \tau,F^N(\bar a) = b
\Leftrightarrow \bar a \in {}^n|N| \wedge b \in N \wedge F^M(\bar a) = b$.
\end{definition}

\begin{definition}
\label{z36}
Let $\cL$ be a logic.

\noindent
0) at$_{\cL}(\tau)$ is the set of atomic formulas $\varphi(\bar x)$
   from \ref{z29}(a); bs$_{\cL}(\tau) = \{\varphi,\neg \varphi:\varphi
\in \text{ at}_{\cL}(\tau)$; we may omit $\cL$ if clear from the context.

\noindent
1) Let $\varphi \vdash_{\cL} \psi$ where $\varphi,\psi \in \cL(\tau)$
mean that $M \models \varphi \Rightarrow M \models \psi$ for any
$\tau$-model (this does not depend on $\tau$ by Definition \ref{z21}).

\noindent
2) We say $\bar\vartheta$ is an $(\cL,\tau_1,\tau_2)$-interpretation
scheme \when :
\begin{enumerate}
\item[$(a)$]  $\cL$ a logic
\sn
\item[$(b)$]  $\tau_1,\tau_2$ vocabularies
\sn
\item[$(c)$]  $\bar\vartheta = \langle \vartheta_{\varphi(\bar
x)}(\bar x):\varphi(\bar x) \in \text{ at}_{\cL}(\tau_2)\rangle$
\sn
\item[$(d)$]  $\vartheta_{\varphi(\bar x)}(\bar x)$ is a formula in
$\cL(\tau_1)$ so $\varphi = P(x_0,\dotsc,x_{n-1}),\bar x = \langle
x_\ell:\ell < n \rangle$ or $\varphi = (x_0 = F(x_1,\dotsc,x_n)),\bar
x = \langle x_\ell:\ell \le  n\rangle$.
\end{enumerate}
\mn
2A) Above we say the $((\cL_1,\tau_1,\tau_2)$-interpretation) scheme
$\bar\vartheta$ is \underline{simple} when
$\vartheta_{x_0 = x_1}$ has the form $(x_0=x_1)$.

\noindent
3) For $\bar\vartheta$ as above we say $N = N_{\bar\vartheta}[M] =
M_1[\bar\vartheta]$ \when \,
\mn
\begin{enumerate}
\item[$(a)$]   $|N| = |M|/\varphi_= (x_0,x_1)$ which means
\sn
\begin{enumerate}
\item[$\bullet$]  $E^M_{\varphi_=} = \{(a,b):M \models \varphi(a,b)\}$
is an equivalence relation on $|M|$ or just some non-empty subset
\sn
\item[$\bullet$]  $|N| = \{a/E^M_{\varphi_=}:a \in M$ and $M \models
  \varphi_=[a,a]\}$ 
\end{enumerate}
\item[$(b)$]   $N \models \varphi[\bar a] \Leftrightarrow M \models
\vartheta_{\varphi(\bar x)}[\bar a]$ for $\varphi(\bar x) \in
\text{ at}(\tau_2)$ and $\bar a \in {}^{\omega >} M$ of length $\ell
g(\bar x)$
\end{enumerate}
\mn
(note: not for every such $\bar\vartheta$ and
$\tau_1$-model $M_1$ is $N_{\bar\vartheta}[M_2]$ well defined, we need
that $\varphi_=(-,-)$ defines an equivalence relation on $|M_1|$, which
is a congruence relation for the $\tau_2$-relations and functions we
define and, of course, the definition of functions 
gives ones, similarly below; but if $\vartheta$ is simple this problem
does not arise).

\noindent
4) We say the logic $\cL$ satisfies full substitution\footnote{We do not use
 full substitution; as for $\kappa$ singular
 $\bbL^0_\kappa,\bbL^1_\kappa$ are not closed under full substitution}
 \when : if $\tau_1,\tau_2$ are vocabularies, $\bar\vartheta
 = \langle \vartheta_{\varphi(\bar x)}(\bar x):\varphi(\bar x) \in
 \text{ at}_{\cL}(\tau_2)\rangle$ is a simple
 $(\cL,\tau_1,\tau_2)$-interpretation scheme, and 
$\psi_2 \in \cL(\tau_2)$, \then \, there is $\psi_1 \in \cL(\tau_1)$ such
 that: if $M_1 = N_{\bar\vartheta}[M_2]$, see below, so $M_\ell$ is a
 $\tau_\ell$-models for $\ell=1,2$ then $M_1 \models \psi_1
 \Leftrightarrow M_2 \models \psi_2$. 

\noindent
5) We say that the logic $\cL$ satisfies substitution 
\when \, we require $\bar\vartheta$ to be simple.
\end{definition}

\begin{definition}
\label{z37}
1) We say $\bar\vartheta$ is a weak $(\cL,\tau_1,\tau_2)$-scheme (but
   $\cL$ is immaterial so can be omitted) \when \,:
\mn
\begin{enumerate}
\item[$(a)$]  $\tau_1,\tau_2$ are vocabularies
\sn
\item[$(b)$]  $\bar\vartheta = \langle \vartheta_{\varphi(\bar
  x)}(\bar x):\varphi(\bar x) \in \at(\tau_2)\rangle$
\sn
\item[$(c)$]  the formulas $\vartheta_{\varphi(\bar x)}$ are atomic
(on conjunction of two atomic for equality) formulas in which we substitute 
some variables by individual constants, moreover:
\sn
\begin{enumerate}
\item[$(\alpha)$]  if $\varphi(\bar x)$ is equal to $(x_0 = x_1)$ so
  $\ell g(\bar x)=2$ then $\vartheta_{\varphi(\bar x)} = (x_0 = x_1
  \wedge P_{\varphi(\bar x)}(x_0,\bar c_{\varphi(\bar x)})$
\sn
\item[$(\beta)$]  if $\varphi(\bar x) = P(\bar x)$, so $P \in \tau_2$
  a predicate \then \,
\sn
\item[${{}}$]  $\bullet_1 \quad \vartheta(\bar x) = Q_{\varphi(\bar
x)}(\bar x,\bar c)$ where
\sn
\item[${{}}$]  $\bullet_2 \quad Q_{\varphi(\bar x)}$ a predicate from $\tau_1$
\sn
\item[${{}}$]  $\bullet_3 \quad \bar c$ a sequence of individual
  constants from $\tau_1$
\sn
\item[${{}}$]  $\bullet_4 \quad$ so arity$_{\tau_2}(P) = \ell g(\bar
  x)$, arity$_{\tau_1}(Q_p) = \ell g(\bar x) + \ell g(\bar c)$
\sn
\item[$(\gamma)$]  if $\varphi(\bar x) = (x_0 = F(x_1,\dotsc,x_n))$
  so $F \in \tau_2$ a function symbol \then \,
\sn
\item[${{}}$]  $\bullet_1 \quad \vartheta_{\varphi(\bar x)}(\bar x)$
  is $(x_0 = H_F(x_1,\dotsc,x_n,\bar c_F))$ or
  $Q_F(x_0,\dotsc,x_n,\bar c_F)$
\sn
\item[${{}}$]  $\bullet_2,\,\bullet_3 \quad$ as above.
\end{enumerate}
\end{enumerate}
\mn
2) We say $\cL$ has weak substitution \when \, as in Definition
\ref{z36}(4), but restricting ourselves to weak schemes.
\end{definition}

\begin{remark}
\label{z38}
1) The meaning of \ref{z37}(1)(c)$(\gamma)$ is
\mn
\begin{enumerate}
\item[$(*)_1$]  $\bullet \quad \vartheta_{\varphi(\bar x)}(\bar x) =
  (x_0 = H_F(x_1,\dotsc,x_n,\bar c_F))$
\sn
\item[${{}}$]  $\bullet \quad H_F \in \tau_1$ a function symbol of
  arity $n + \ell g(\bar c_F)$
\sn
\item[${{}}$]  $\bullet \quad \bar c_F$ a sequence of individual
  constants from
\newline
\underline{or}
\sn
\item[$(*)_2$]  similarly using $Q'_F(x_0,\dotsc,x_n,\bar c_F)$.
\end{enumerate}
\mn
2) Why in \ref{z37}(1)(c) the ``moreover", i.e. why is the second
version stronger?  As we demand the $\bar c_\varphi$'s to be as in the
end.  Does not matter whether we ask it or not.  We could further
demand $\bar c_\varphi = \langle c_* \rangle$ for all $\bar\varphi \in
\at(\tau_2)$. 

\noindent
3) In ``weak $(\cL,\tau_1,\tau_2)$-scheme" (hence in ``weak
substitution") we may use first order formulas (instead of atomic),
   i.e. that is our results will not be affected by this change in the
   definition.

\noindent
4) Also, in ``weak$(\cL,\tau_1,\tau_2)$-scheme" we may add $\bar
c_\varphi = \langle c_* \rangle$, for one $c_*$, and/or demand
$\bar\vartheta$ is simple.  For the later change in the proof of
\ref{a28}, we have to say that ``\wilog \, all $\|M_{\alpha,n}\|$ have
cardinality $\ge \mu_1$ and in fact $= \mu_1$", for this we need claim
version of LST, see so have to add the assumption ``$\cL$ satisfies
suitable version of LST.
\end{remark}

\begin{definition}
\label{z39}
1) For logics $\cL_1,\cL_2$ let $\cL_1 \subseteq \cL_2$ mean that:
   $\cL_1(\tau) \subseteq \cL_2(\tau)$ for any vocabulary $\tau$ and
   $M \models_{\cL_1} \varphi$ is equivalent to $M \models_{\cL_2}
   \varphi$ when $\varphi \in \cL(\tau)$ and $M$ is a $\tau$-model.

\noindent
2) For logics $\cL_1,\cL_2$ let $\cL_1 \le \cL_2$ or $\cL_2$ is
stronger than $\cL_2$ means that for every $\tau$ and $\varphi \in
\cL_1(\tau)$ there is $\psi \in \cL_2(\tau)$ such that $\varphi,\psi$
are equivalent, i.e. $M \models_{\cL_1} \varphi \Leftrightarrow M
\models_{\cL_2} \psi$ for any $\tau$-model $M$.

\noindent
3) We say that the logics $\cL_1,\cL_2$ are equivalent,
$\cL_1 \equiv \cL_2$ when $\cL_1 \le \cL_2$ and $\cL_2 \le \cL_1$.
\end{definition}

\begin{definition}
\label{z42}
1) The logic $\bbL_{\kappa,\theta}$ for $\kappa 
\ge \theta \ge \aleph_0$ is defined like first order logic, but 
$\bbL_{\kappa,\theta}(\tau)$ is the closure of the set of
   atomic formulas under $\neg \varphi,\bigwedge\limits_{i < \alpha}
   \varphi_i$ where $\alpha < \kappa$ and $(\forall
   x_0,\dotsc,x_i,\ldots)_{i < \alpha} \varphi$ where $\alpha <
   \theta$; a sentence is a formula with no free variables and
   satisfaction is defined naturally.

\noindent
1A) First order logic is $\bbL = \bbL_{\aleph_0,\aleph_0}$.

\noindent
2) Let $\bbL_{\kappa,\theta,\gamma}(\tau)$ be the set of formulas
   $\varphi(\bar x) \in \bbL_{\kappa,\theta}(\tau)$ of depth $\le
   \gamma$; similarly in the other cases.

\noindent
3) Let $\bbL^0_{\le \kappa,\le \theta} = \bbL_{\kappa^+,\theta^+},
\bbL^0_{< \kappa,< \theta} = \cup\{\bbL^0_{\le \kappa_1,\le
\theta_1}:\kappa_1 < \kappa,\theta_1 < \theta\}$; similarly $\bbL^0_{\le
\kappa,\le \theta,< \gamma}$, etc.

\noindent
4) Let $\bbL^{-1}_\kappa = \bbL^0_{< \kappa,\aleph_0},
\bbL^0_\kappa = \bbL^0_{< \kappa,< \kappa}$.

\noindent
5) We define the logic $\bbL_{\infty,\theta,\gamma}$ by induction on
   the ordinal $\gamma$ such that $\bbL_{\infty,\theta,\gamma}(\tau)$
   is a set of cardinality $\le \beth_\gamma(|\tau| + \aleph_0)$,
   increasing with $\gamma$ as follows:
\mn
\begin{enumerate}
\item[$\bullet$]  for $\gamma =0$, the set of basic sentences
\sn
\item[$\bullet$]  for $\gamma$ limit, it is
$\cup\{\bbL_{\infty,\theta,\beta}(\tau):\beta < \gamma\}$
\sn
\item[$\bullet$]  for $\gamma = \beta +1$ it is the set of sentences
of the form $(\exists \bar x) \bigwedge\limits_{i < \alpha}
\varphi_i(\bar x)$ or its negation where $\ell g(\bar x)
<\theta,\varphi_i(\bar x)$ a formula from
$\bbL_{\infty,\theta,\beta}(\tau)$ with $\alpha < \beth_\beta(|\tau| +
\aleph_0)$.
\end{enumerate}
\end{definition}

\begin{observation}
\label{z45}
If $\kappa$ is strong limit singular, \then \, there is no logic $\cL$ such
that:
\mn
\begin{enumerate}
\item[$(a)$]  $\bbL^{-1}_\kappa \le \cL \le \bbL^0_\kappa$
\sn
\item[$(b)$]  $\cL$ satisfies full substitution.
\end{enumerate}
\end{observation}

\noindent
Note that
\begin{claim}
\label{z47} 
1) In Definition \ref{z42} for $\kappa$ singular, $\bbL_{\kappa^+,\theta} \le
   \bbL_{\kappa,\theta} \subseteq \bbL_{\kappa^+,\theta}$, i.e. every
   sentence in $\bbL_{\kappa^+,\theta}(\tau)$ is equivalent to one in
   $\bbL_{\kappa,\theta}$ (and $\bbL_{\kappa,\theta}(\tau) \subseteq
   \bbL_{\kappa^+,\theta}(\tau))$, so it seems pointless to allow
   $\kappa$ to be singular as was originally done.  Still this is the
   tradition and we use it.

\noindent
2) If $\kappa$ is a regular cardinal \then \, the logics
$\bbL^0_{\kappa,\theta}$ and $\bbL^0_{< \kappa,\theta}$ are equal.

\noindent
3) If $\kappa = \beth_\kappa$ then $\bbL^0_{< \kappa,<\kappa,<\kappa}
   = \bbL^0_{< \kappa,<\kappa}$.
\end{claim}
\newpage

\section {The logic $\bbL^1_\kappa$}

Our main definition is
\begin{definition}
\label{a8}
For a vocabulary $\tau,\tau$-models $M_1,M_2$, set $\Gamma$ 
formulas in the vocabulary $\tau$ in any logic (each with
finitely many free variables if not said otherwise (see \ref{a10}(4)),
cardinal $\theta$ and ordinal $\alpha$ we define a game $\Game =
\Game_{\Gamma,\theta,\alpha}[M_1,M_2]$ 
as follows, and using $(M_1,\bar b_1),(M_2,b_2)$ with their natural meaning
when Dom$(\bar b_1) = \text{ Dom}(\bar b_2)$.
\mn
\begin{enumerate}
\item[$(A)$]  The moves are indexed by $n < \omega$ (but every actual play is
finite), just before the $n$-th move we have a state $\bold
s_n = (A^1_n,A^2_n,h^1_n,h^2_n,g_n,\beta_n,n)$
\sn
\item[$(B)$]  $\bold s = (A^1,A^2,h^1,h^2,g,\beta,n) = (A^1_{\bold
s},A^2_{\bold s},h^1_{\bold s},h^2_{\bold s},g_{\bold s},
\beta_{\bold s},n_{\bold s})$ is a state (or $n$-state or
$(\theta,n)$-state or $(\theta,< \omega)$-state) \when :
\begin{enumerate}
\item[$(a)$]   $A^\ell \in [M_\ell]^{\le \theta}$ for $\ell=1,2$
\sn
\item[$(b)$]   $\beta \le \alpha$ so an ordinal
\sn
\item[$(c)$]  $h^\ell$ is a function from $A^\ell$ into $\omega$
\sn
\item[$(d)$]   $g$ is a partial one-to-one function from $M_1$ to $M_2$ and
let $g^1_s = g^1 = g_s =g$ and let $g^2_s = g^2 = (g^1_s)^{-1}$,
\sn
\item[$(e)$]  Dom$(g^\ell) \subseteq A^\ell$ for $\ell=1,2$
\sn
\item[$(f)$]  $g$ preserves satisfaction of the formulas in $\Gamma$ and
their negations, i.e. for $\varphi(\bar x) \in \Gamma$ and $\bar a \in
{}^{\ell g(\bar x)}\text{Dom}(g)$ we have $M_1 \models \varphi[\bar
a] \Leftrightarrow M_2 \models \varphi[g(\bar a)]$
\sn
\item[$(g)$]  if $a \in \text{ Dom}(g^\ell)$ then $h^\ell(a) < n$
\end{enumerate}
\item[$(C)$]  we define the state $\bold s = \bold s_0 = \bold
s^0_\alpha$ by letting $n_{\bold s} = 0,A^1_{\bold s} = \emptyset = 
A^1_{\bold s},\beta_{\bold s} = \alpha,h^1_{\bold s} = \emptyset 
= h^2_{\bold s},g_s = \emptyset$; so really $\bold s$ depend only on $\alpha$ 
(but in general, this may not be a state for our game as possibly for 
some sentence $\psi \in \Gamma$ we have
$M_1 \models \psi \Leftrightarrow M_2 \models \neg \psi$)
\sn
\item[$(D)$]  we say that a state $\bold t$ extends a state $\bold s$
when $A^\ell_{\bold s} \subseteq A^\ell_{\bold t},h^\ell_{\bold s}
\subseteq h^\ell_{\bold t}$ for $\ell=1,2$ and $g_{\bold s} \subseteq
g_{\bold t},\beta_{\bold s} > \beta_{\bold t},n_{\bold s} < n_{\bold
t}$; we say $\bold t$ is a
successor of $\bold s$ if in addition $n_{\bold t} = n_{\bold s} +1$
\sn
\item[$(E)$]  in the $n$-th move

\underline{the anti-isomorphism} player (AIS) chooses
 $(\beta_{n+1},\iota_n,A'_n)$ such that:

$\iota_n \in \{1,2\},
\beta_{n+1} < \beta_n$ and $A^{\iota_n}_n \subseteq A'_n 
\in [M_{\iota_n}]^{\le \theta}$, 

\underline{the isomorphism} player (ISO) chooses a state $\bold s_{n+1}$
such that
\begin{enumerate}
\item[$\bullet$]   $\bold s_{n+1}$ is a successor of $\bold s_n$
\sn
\item[$\bullet$]  $A^{\iota_n}_{\bold s_{n+1}} = A'_n$
\sn
\item[$\bullet$]  $A^{3-\iota_n}_{\bold s_{n+1}} = 
A^{3-\iota_n}_{\bold s_n} \cup \Dom(g^{3-\iota_n}_{\bold s_{n+1}})$
\sn
\item[$\bullet$]  if $a \in A'_n \backslash A^{\iota_n}_{\bold s_n}$ then
$h^{\iota_n}_{\bold s_{n+1}}(a) \ge n+1$
\sn
\item[$\bullet$]  $\Dom(g^{\iota_n}_{\bold s_{n+1}}) = \{a \in
A^{\iota_n}_{\bold s_n}:h^{\iota_n}_{\bold s_n}(a) < n+1\}$ so it includes
$\Dom(g^{\iota_n}_{\bold s_n})$ 
\sn
\item[$\bullet$]  $\beta_{\bold s_{n+1}} = \beta_{n+1}$.
\end{enumerate}
\item[$(F)$]  $\bullet \quad$ the play ends when one of the player has no 
legal moves (always occur

\hskip25pt as $\beta_n < \beta_{n-1}$) and then this player loses, this
may occur for $n=0$
\sn
\item[${{}}$]  $\bullet \quad$ for $\alpha = 0$ we stipulate that ISO
wins iff $\bold s^0_\alpha$ is a state.
\end{enumerate}
\end{definition}  

\begin{discussion}
\label{a9}
1) This is a parallel to EF-games.  Note that we like on the one hand
   the game to have $\le \omega$ moves, really each play has $< \omega$
   moves and deal with sets of cardinality $\le \theta$ and 
on the other hand we do not like well ordering to be definable, 
i.e. allow $M_1$ to be well ordered  while $M_2$ to be non-well
 ordered but still the ISO player wins.  We do 
this by ``rescheduling our debts", i.e. using the $h_n$'s.
\end{discussion}

\begin{definition}
\label{a10}
1) Let $\cE^{0,\tau}_{\Gamma,\theta,\alpha}$ be the class
$\{(M_1,M_2):M_1,M_2$ are $\tau$-models and in the game
$\Game_{\Gamma,\theta,\alpha}[M_1,M_2]$ the ISO player has a winning
   strategy$\}$ where $\Gamma$ is a set of formulas in the vocabulary
   $\tau$, each with finitely many free variables.

\noindent
2) $\cE^{1,\tau}_{\Gamma,\theta,\alpha}$ is the closure of 
$\cE^{0,\tau}_{\Gamma,\theta,\alpha}$ to an equivalence relation (on the
class of $\tau$-models).

\noindent
3) Above we may replace $\Gamma$ by qf$(\tau)$ which means $\Gamma =$
   the set at$(\tau)$ or bs$(\tau)$ formulas in the vocabulary $\tau$.

\noindent
4) Above if we omit $\tau$ we mean $\tau = \tau_\Gamma$ and if we omit
   $\Gamma$ we mean bs$(\tau)$.  Abusing
   notation we may say $M_1,M_2$ are
   $\cE^{0,\tau}_{\Gamma,\theta,\alpha}$-equivalent. 
\end{definition}

\begin{fact}
\label{a12}
Assume $\Game_{\Gamma,\theta,\alpha}[M_1,M_2]$ is well defined and
$M_1,M_2$ are $\tau$-models.

\noindent
1) The game $\Game_{\Gamma,\theta,\alpha}[M_1,M_2]$ is a determined
game and is 
without memory, i.e. during a play, being a winning situation does not
depend on the history, just on the current state, also only $M_\ell 
\rest \tau_\Gamma$ are relevant.

\noindent
2a) The relation $\cE^0_{\Gamma,\theta,\alpha}$ holds for
$(M_1,M_2)$ when $M_1,M_2$ are isomorphic $\tau$-models.

\noindent
2b) If $M_1 \cong M'_1,M_2 \cong M'_2$ then
$M_1 \cE^0_{\Gamma,\theta,\alpha} M_2 \Leftrightarrow M'_1
\cE^0_{\Gamma,\theta,\alpha} M'_2$.

\noindent
2c) $\cE^0_{\Gamma,\theta,\alpha}$ is reflexive and symmetric.

\noindent
3) The relation $\cE^1_{\Gamma,\theta,\alpha}$ is an  
equivalence relation on the class of $\tau$-models.

\noindent
4) If $\alpha$ is a limit ordinal \then \, $M_1
\cE^0_{\Gamma,\theta,\alpha} M_2$ iff $[\beta < \alpha
\Rightarrow M_1 \cE^0_{\Gamma,\theta,\beta} M_2]$.

\noindent
5) $\cE^1_{\Gamma,\theta,\alpha}$ has $\le \beth_{\alpha +1}(|\Gamma|+\theta)$
   equivalence classes.

\noindent
6) If $\tau^1 \subseteq \tau^2,\Gamma_1 \subseteq \Gamma_2,
\theta_1 \le \theta_2$ and $\alpha_1 \le \alpha_2$ and $M_1
\cE^{0,\tau^2}_{\Gamma_2,\theta_2,\alpha_2} M_2$ \then \, $M_1
\cE^{0,\tau^1}_{\Gamma_1,\theta_1,\alpha_1} M_2$.
\end{fact}

\begin{PROOF}{\ref{a12}}
1) Obvious.

\noindent
2a) Let $g_*$ be an isomorphism from $M_1$ onto $M_2$.  Now a winning
   strategy for the player ISO in
   $\Game_{\Gamma,\theta,\alpha}[M_1,M_2]$ is to preserve ``$g_{\bold
   s_n} \subseteq g_*$".

\noindent
2b) Should be clear.

\noindent
2c) Let us check.

\underline{Reflexivity}:

Follows from part (2a) as $M \cong M$.

\underline{Symmetry}:

Reading the definition carefully it should be clear.

\noindent
3),4) Easy, too.

\noindent
5)  We prove by induction on $\alpha$ that there is an equivalence
relation $\cE^2_{\Gamma,\theta,\alpha}$ on the class of
$\tau(\Gamma)$-models such that $M_1 \cE^2_{\Gamma,\theta,\alpha} M_2
\Rightarrow M_1 \cE^0_{\Gamma,\theta,\alpha} M_2$ and
$\cE^2_{\Gamma,\theta,\alpha}$ has $\le \beth_{\alpha +1}(|\Gamma| +
\theta)$ equivalence classes.  For $\alpha = 0$, recall
 clause (F) of Definition \ref{a8} so easily $M_1 \cE^0_{\Gamma,\theta,\alpha}
   M_2$ means that $\Gamma_\ell = \{\varphi[\bar a]:M_\ell \models
   \varphi$ and $\varphi \in \Gamma,\varphi$ is a sentence$\}$ does
   not depend on $\ell$; clearly $\cE^0_{\Gamma,\theta,\alpha}$ 
is an equivalence relation with $\le 2^{|\Gamma|}$ equivalence classes
and $\cE^1_{\Gamma,\theta,\alpha} = \cE^0_{\Gamma,\theta,\alpha}$; let
$\cE^2_{\Gamma,\theta,\alpha} = \cE^0_{\Gamma,\theta,\alpha}$.
  For $\alpha$ a limit ordinal use part (4) and choose
$\cE^2_{\Gamma,\theta,\alpha} =
\cap\{\cE^2_{\Gamma,\theta,\beta}:\beta < \alpha\}$ 
recalling $\beth_\alpha(|\Gamma| + \theta) =
   \Sigma\{\beth_\beta(|\Gamma| + \theta):\beta < \alpha\}$.

Lastly, for $\alpha = \beta +1$ use the induction hypothesis and
$\beth_\alpha(|\Gamma| + \theta) = 2^{\beth_\beta(|\Gamma| +
\theta)}$.  Alternatively use $\bbL_{\theta^+,\theta^+,\alpha}$,
i.e. use the proof of \ref{a14}, when we replace
$\beth_\alpha(|\Gamma| + \theta)$ by $\beth_{(|\Gamma| + \theta)^+}$.

\noindent
6) Easy.
\end{PROOF}

\begin{definition}
\label{a13}
We define the logic $\bbL^1_{\le \theta}$ as follows:
a sentence $\psi \in \bbL_{\le \theta}(\tau)$ \underline{iff} the 
sentence is defined using (or by) a triple
$(\text{qf}(\tau_1),\theta,\alpha)$ which means:
$\tau_1$ a sub-vocabulary of $\tau$ of
cardinality $\le \theta$ and $\alpha < \theta^+$ and for 
some sequence $\langle M_\alpha:\alpha < \alpha(*)\rangle$
of $\tau_1$-models of length $\alpha(*) \le \beth_{\alpha +1}(\theta)$ 
we have: $M \models \psi$ \underline{iff} $M$ is
$\cE^1_{\text{qf}(\tau_1),\theta,\alpha}$-equivalent to $M_\alpha$ for
some $\alpha < \alpha(*)$. 

\noindent
2) Let $\bbL^1_\kappa = \cup\{\bbL^1_{\le \theta}:\theta < \kappa\}$
   so $\bbL^1_{\theta^+} = \bbL^1_{\le \theta}$.
\end{definition}

\begin{remark}
\label{a13c}
1) The present definition of the logic (\ref{a13}) is interesting mainly 
for $\kappa$ strong limit such that 
$\kappa = \beth_\kappa$ and it seems to me that it makes  
the presentation transparent.  Note that $\bbL^1_{\le \theta}$ is
similar to the set of formulas of quantifier depth $< \theta^+$.  

\noindent
2) Why? Just note that for vocabulary $\tau$ of cardinality $\le
\theta,\bbL^1_{\le \theta}(\tau)$ has cardinality $\beth_{\theta^+}$, and
it helps to have arbitrary Boolean combinations of formulas of a fix
quantifier depth $[< \theta^+]$.  
Note that if $\varphi(\bar x) \in \bbL^1_\kappa(\tau)$ has
infinitely many free variables, we cannot ``close" it to a sentence.

\noindent
3) We may instead define $\bbL^{1,*}_{\le \theta}$ by: $\psi \in
   \bbL^1_{\le \theta}(\tau)$ iff for some $\tau' \subseteq \tau$ of
   cardinality $\le \theta$ and some sequence $\langle M_i:i <
i(*)\rangle$ of $\tau'$-models of cardinality $\le 2^\theta$ and
   some $\gamma <\theta^+$ we have: a $\tau$-model $M$ satisfies
$\psi$ iff $M \cE^1_{\text{qf}(\tau'),\theta,\gamma} M_i$ for some $i$.  
\end{remark}

\begin{claim}
\label{a13m}
1) $\bbL^1_{\le \theta}$ is a nice logic, see Definition \ref{z21}.

\noindent
2) The logic $\bbL^1_{\le \theta}$ has full substitution.

\noindent
3) $\bbL^1_{\le \theta}(\tau)$ has cardinality $\le |\tau|^\theta +
\beth_{\theta^+}$ for any vocabulary $\tau$.

\noindent
4) If $\kappa = \beth_\kappa$ \then \, $\bbL^1_\kappa$ is a nice logic,
$\bbL^1_\kappa(\tau)$ has cardinality $\kappa$ whenever $\tau$
 is a vocabulary of cardinality $< \kappa$ and $\bbL^1_\kappa(\tau)$ has
 cardinality $|\tau|^{< \kappa} + \kappa$ for any vocabulary $\tau$.

\noindent
5) In part (4), if $\kappa$ is regular \then \, the logic $\bbL^1_\kappa$
has full substitution.

\noindent
6) If $\theta,\alpha < \kappa$ and $\tau$ is a vocabulary of
   cardinality $< \kappa$ \then \, for any set $\cU$ of
$\cE^1_{\text{\rm qf}(\tau),\theta,\alpha}$-equivalence classes for some
$\psi \in \bbL^1_\kappa(\tau)$ we have 
$\{M:M/\cE^1_{\text{\rm qf}(\tau),\theta,\alpha} \in \cU\} = \{M:M$ a
   $\tau$-model of $\psi\}$.  

\end{claim}

\begin{PROOF}{\ref{a13m}}
1),2)  Just check definition \ref{z21}, \ref{z29}, \ref{z36}(5) but still
we elaborate the relatively more substantial \ref{z29}, \ref{z36}(5).

\noindent
\underline{Clause $(a)$}:  \underline{Atomic formulas}

Just note that if $M_1 \cE^0_{\text{qf}(\tau),\theta,\alpha}
M_2$ for $\alpha= 0$ (any $\theta$) \then \, $M_1,M_2$ satisfies the
same $\tau$-atomic sentences.

\noindent
\underline{Clause $(b)(\alpha)$}:  \underline{Conjunction, 
(similarly Disjunction)}

Just note that if $\varphi_\ell$ is defined using 
$(\text{qf}(\tau_\ell),\theta_\ell,\alpha_\ell)$ and $\tau = \tau_1
\cup \tau_2,\theta = \text{ max}\{\theta_1,\theta_2\},\alpha = \text{
max}\{\alpha_1,\alpha_2\}$ then the equivalence relation
$\cE^1_{\text{qf}(\tau),\theta,\alpha}$ refine the equivalence
relation $\cE^1_{\text{qf}(\tau_\ell),\theta_\ell,\alpha_\ell}$ for
$\ell=1,2$.

\noindent
\underline{Clause $(b)(\beta)$}:  \underline{Existential Quantifier}

Assume $\varphi \in \bbL^1_\kappa(\tau \cup \{c\})$ is defined by
the triple $(\qf(\tau \cup \{c\}),\theta,\alpha)$ and $\varphi(x)$ is
the corresponding formula $\bbL^1_\kappa(\tau)$ and $\exists x
\varphi(x)$ is the naturally defined sentence.

Now if $M_1,M_2$ are $\cE^0_{\qf(\tau),\theta,\alpha
+1}$-equivalent $\tau$-models then $M_1 \models \exists x \varphi(x)$
iff $M_2 \models (\exists x)\varphi(x)$ by the definition of the game.
Hence this holds for ``$M_1,M_2$ are
$\cE^1_{\qf(\tau),\theta,\alpha +1}$-equivalent $\tau$-models".

\noindent
\underline{Clause $(b)(\gamma)$}:  \underline{Negation}

Obvious by the definition because for any vocabulary $\tau$ of
cardinality $\le \theta$ the equivalence relation
$\cE^1_{\text{qf}(\tau),\theta,\alpha}$ has $\le \beth_{\alpha
+1}(\theta)$ equivalence classes by \ref{a12}(5).

\noindent
\underline{Clause $(c)$}:  \underline{Restricting to a unary predicate $P$}

Easily if $\psi$ is defined using $(\qf(\tau),\theta,\alpha)$ and
\wilog \, $P \in \tau$ then so is $\psi \rest P$.

\noindent
\underline{Clause $(d)$}:  \underline{We prove more: Full substitution}:  

Assume we are given vocabularies $\tau_1,\tau_2$ and consider substituting
$\varphi_P(x_0,\dotsc,x_{\text{arity}_{\tau_2}(P)-1}) \in \bbL^1_{\le
\theta}(\tau_0)$ for $P \in \tau_2$ treating 
$F(x_0,\dotsc,x_{\text{arity}(F)-1}) = x_{\text{arity}(F)}$ as an 
$(\text{arity}_{\tau_2(F)+1})$-place predicate.
Let $\varphi_P$ be defined by
$(\text{qf}(\tau^{+\text{arity}(P)}),\theta_P,\alpha_P)$ and
let $\alpha_0 = \sup\{\alpha_P:P \in \tau_2\}$.

We are given $\psi \in \bbL^1_{\le \theta}(\tau_1)$ and we shall find $\psi'
\in \bbL^1_{\le \theta}(\tau_2)$ which says that if we substitute
$\varphi_P(x_0,\dotsc,x_{\text{arity}_{\tau_2}(P)})$ 
instead $P(x_0,\ldots)$ in $\psi$ for every $P \in
\tau_1$, we get (up to equivalence) $\psi'$.  Let $\psi$ be defined by
$(\text{qf}(\tau_1),\theta_1,\alpha_1)$.  Let $\alpha_2 = \alpha_1 +
\alpha_0$ and easily there is $\psi'$ as required defined by
$(\text{qf}(\tau_2),\theta_1,\alpha_2)$. 
\end{PROOF}

\begin{claim}
\label{a14}
1) Let $\kappa = \beth_{\theta^+}$.

\noindent
$\bbL^1_{\le \theta} \le \bbL^0_\kappa$, i.e. every 
formula of $\bbL^1_{\le \theta}$ is equivalent to, hence
 can be looked at, as a formula of $\bbL^0_\kappa$.

\noindent
2) $\bbL^{-1}_{\le \theta} \le \bbL^1_{\le \theta}$.
\end{claim}

\begin{remark}
\label{a15}
For many purposes we identify them, i.e. say $\bbL^1_{\le \theta}
\subseteq \bbL^0_\kappa$; this gives
another reasonable version of subformulas.
\end{remark}

\begin{PROOF}{\ref{a14}}
1) We first prove:
\mn
\begin{enumerate}
\item[$\boxplus_1$]  if $\bold s$ is a state in the game
$\Game_{\qf(\tau),\theta,\alpha}[M_1,M_2]$ and $\beta =
\beta_{\bold s},\tau = \tau(M_\ell),\bar a_1 = \langle
a^1_\varepsilon:\varepsilon < \varepsilon(*)\rangle \in
{}^{\theta^+>}(M_1)$ list $\Dom(g_{\bold s})$ and $\bar a_2 = \langle
a^2_\varepsilon:\varepsilon < \varepsilon(*)\rangle$ where
$a^2_\varepsilon = g_{\bold s}(a^1_\varepsilon)$ and $M_2 \models
\varphi[\bar a_1] \Leftrightarrow M_2 \models \varphi[\bar a_2]$ for
every $\varphi = \varphi(\langle x_\varepsilon:\varepsilon <
\varepsilon(*)\rangle) \in L^*_\beta := \bigcup\limits_{\zeta <
\beta} \bbL_{(\beth_\zeta(\theta + |\tau|))^+,\theta^+}(\tau)$, or just
$\varphi \in L^*_\beta =\bbL_{\infty,\theta^+,\beta}(\tau)$, see \ref{z42}(5),
 \then \, $\bold s$ is a winning state for the player ISO in the game.
\end{enumerate}
\mn
[Why?  We prove this by induction on $\beta$.  First, for $\beta = 0$
this is trivial as $\bold s$ is a state.  Second, for $\beta$ limit,
any choice of the AIS player includes an ordinal $\gamma =
\beta_{n_{\bold s}+1} < \beta$, so
the ISO player may ``pretend" that the given state $\bold s$ has
$\beta_{\bold s} = \gamma +1$ and use the induction hypothesis.
Third, if $\beta = \gamma +1$, let the AIS player make his choice
$(\beta_{n_{\bold s}+1},\iota,A)$.  
Now ISO has to extend $g^\iota_{\bold s}$ adding some
$\le \theta$ elements of $M_\iota$ to its domain, the elements in $\{a
\in \Dom(g^\iota_{\bold s}):h^\iota_{\bold s}(a) < n_{\bold s} +1\}
\cup A$, let $\bar b_\iota$ list them.  
Let $\bar x = \langle x_\varepsilon:\varepsilon <
\varepsilon(*)\rangle$ and let $\bar y = \langle y_\varepsilon:\varepsilon
< \ell g(\bar b_\iota)\rangle$ and define $\varphi_*(\bar y,\bar x) =
\bigwedge\{\varphi(\bar y,\bar x) \in L^*_\gamma:M_\iota \models
\varphi[\bar b,\bar a_\iota]\}$.

So $M_\iota \models \varphi_*[\bar b,\bar a_\iota]$ hence $M_\iota
\models (\exists \bar y)\varphi_*(\bar y,\bar a_\iota)$ hence by the
assumption of $\boxplus_1$ we have $M_{3-\iota} \models ``(\exists \bar
y)\varphi_*(\bar y,\bar a_{3-\iota})"$ hence there is $\bar b_{3-\iota} \in
{}^{\ell g(\bar b_\iota)}(M_{3-\iota})$ such that $M_{3-\iota} \models
\varphi_*[\bar b_{3-\iota},\bar a_{3-\iota}]$.

Now the player ISO can make its move getting the state $\bold t$, a
successor of $\bold s$ such
that $g_{\bold t} = g_{\bold s} \cup \{\langle
b^1_\varepsilon,b^2_\varepsilon\rangle:\varepsilon < \ell g(\bar
y)\rangle,g_{\bold t}$ maps $\bar b_1$ to $\bar b_2$ and
$n_{\bold t} = n_{\bold s}+1$ and $\beta_{\bold t} = \gamma$;
 clearly possible.  As the ordinal $\beta_{n_{\bold s}+1}$ chosen
by the AIS is $\le \gamma$ also $\bold t$ is as required in
$\boxplus_1$ but so, by the induction hypothesis also it is a
winning state.]

Let
\mn
\begin{enumerate}
\item[$\boxplus_2$]  $\cE_{\tau,\theta,\alpha} = \{(M_1,M_2):M_1,M_2$
are $\tau$-models satisfying the same sentences from $L^*_\alpha\}$.
\end{enumerate}
\mn
Now we shall prove $\boxplus_3$ and by $\boxplus_3(d)$ we are done.
\mn
\begin{enumerate}
\item[$\boxplus_3$]  $(a) \quad \cE_{\tau,\theta,\alpha}$ is an
equivalence relation with $\le (\beth_{\alpha +1}(\theta + |\tau|))$
equivalence 

\hskip25pt classes
\sn
\item[${{}}$]  $(b) \quad \cE_{\tau,\theta,\alpha}$
is an equivalence relation including $\{(M_1,M_2):M_1,M_2$ are

\hskip25pt $\tau$-models such that the player ISO wins in
$\Game_{\text{qf}(\tau),\theta,\alpha}[M_1,M_2]\}$
\sn
\item[${{}}$]  $(c) \quad \cE_{\tau,\alpha,\theta}$ refines the
equivalence relation $\cE^1_{\text{qf}(\tau),\theta,\alpha}$ 
\sn
\item[${{}}$]  $(d) \quad \cE^1_{\text{qf}(\tau),\theta,\alpha}$ has $\le
\beth_{\alpha +1}(\theta + |\tau|)$-equivalence classes.
\end{enumerate}
\mn
[Why?  Clause (a) follows from the number of such sentences being $\le
\beth_\alpha(\theta + |\tau|)$, clause (b) follows by $\boxplus_1$, clause (c)
follows from (b) and the definitions of
$\cE^0_{\text{qf}(\tau),\theta,\alpha}$ and
$\cE^1_{\text{qf}(\tau),\theta,\alpha}$.  Lastly, clause (d) 
follows from clauses (b) + (c).

\noindent
2) Easy.
\end{PROOF}

\begin{remark}
\label{a17}
1) The proof above may seem wasteful
\underline{but} for our purposes this is immaterial.

\noindent
2) We can do better as follows, i.e. another proof of \ref{a14} 
runs as follows:

\noindent
By induction on the formula.  Without loss of
generality we deal with a formula $\Theta_{N,\bar d,
\theta,\alpha,\Gamma}(\bar x)$ by induction on $\beta \le \alpha$ we prove:
\mn
\begin{enumerate}
\item[$\boxplus_1$]  if $n \in [-1,\omega)$ and $\gamma(1),\gamma(2) 
< \theta^+,\bar b \in {}^{\gamma(2)} N,\bold h_\ell:\gamma(\ell) 
\rightarrow \omega$ letting $\bar x = \langle x_i:i <
\gamma(2)\rangle;u_\ell \subseteq \gamma(\ell)$ for $\ell = 1,2,\bold
g$ a one-to-one mapping from $u_1$ onto $u_2$  (so $n= -1 \Rightarrow
\gamma(1),\gamma(2) = 0$) \then \, for some formula $\varphi(\bar z,\bar x) = 
\varphi^0_{N,\bar d,\theta,\alpha,\Gamma,\bold h_1,\bold h_2,g,\beta}(\bar x) 
\in \bbL^0_\kappa$ we have:  for any model $M$
of the vocabulary of $N$ and $\bar c \in {}^{(\ell g(\bar d))}M$ 
and $\bar a \in {}^{\ell g(\bar b)}M$ letting
\begin{enumerate}
\item[$\bullet$]  $A^1_n = \text{ Rang}(\bar c),A^2_n = \text{ Rang}(\bar d)$
\sn
\item[$\bullet$]  $g_n$ has domain $\{c_i:i \in u\}$ and $g_n(c_i) =
d_{\bold g(i)}$ for $i \in u$
\sn
\item[$\bullet$]  $h^\ell_n(c_i) = d_{\bold h_\ell}(i)$
\end{enumerate}
the following are equivalent:
\begin{enumerate}
\item[$(i)$]  $M \models \varphi[\bar c,\bar a]$
\sn
\item[$(ii)$]   $(A_n,B_n,h^1_n,h^2_n,g_n,\beta)$ is a winning position for the
equivalence player in the game $\Game_{N,\theta,\alpha,\Gamma}[(M,\bar
c),(N,\bar d)]$
\end{enumerate}
\item[$\boxplus_2$]   similarly a formula
$\varphi^1_{N,\bar c,\theta,\alpha,\Gamma,h,h',\beta}(\bar z,\bar x)$
expressing situation after the AIS player moves.
\end{enumerate}
\mn
The proof is straightforward.
\end{remark}

\begin{cc}
\label{a18}
We have $M_n \equiv_{\bbL^1_\theta} M_\omega$ for $n < \omega$ and
even $M_n \models \psi[\bar a] \Leftrightarrow M_\omega \models \psi[\bar a]$
\when :
\mn
\begin{enumerate}
\item[$(a)$]  $\psi(\bar z) \in \bbL^1_{\le \theta}(\tau)$ a formula
\sn
\item[$(b)$]  $M_n \prec_{\bbL_{\partial^+,\theta^+}} M_{n+1}$ where
$\partial = \beth_{\theta^+}$
\sn
\item[$(c)$]  $M_\omega := \bigcup\limits_{n < \omega} M_n$
\sn
\item[$(d)$]  $\bar a \in {}^{\ell g(\bar z)}(M_0)$
\sn
\item[$(e)$]  $\tau = \tau(M_n)$ for $n < \omega$.
\end{enumerate}
\end{cc}

\begin{remark}
In fact the logic $\bbL^0_{\partial^+,\theta^+,\theta^+}$ suffice, see
Definition \ref{z42}(2).
\end{remark}

\begin{PROOF}{\ref{a18}}
Without loss of generality $|\tau| \le \theta$.  We shall prove by
induction on $\beta < \theta^+$ the following
\mn
\begin{enumerate}
\item[$\boxplus_\beta$]  $\bold s$ is a winning state for the ISO
player in the game $\Game_{\qf(\tau),\theta,\alpha}[M_n,M_\omega]$ 
\when \, for some $k$ we have
\begin{enumerate}
\item[$(a)$]  $\bold s$ is a state in this game
\sn
\item[$(b)$]  $n \le k < \omega$
\sn
\item[$(c)$]  $\beta_{\bold s} = \beta$
\sn
\item[$(d)$]  Rang$(g^2_{\bold s}) \subseteq M_k$
\sn
\item[$(e)$]  if $b \in A^2_{\bold s}$ then $b \in M_{h^2_{\bold s}(b)}$
\sn
\item[$(f)$]  letting $\bar a^1$ list Dom$(g_{\bold s})$ and $\bar a^2$ be
$\langle g_{\bold s}(a^1_i):i < \ell g(\bar a^1)\rangle$ and $\lambda_\beta =
\beth_\beta(\theta)$ for every $\varphi(\bar x) \in
\bbL_{\lambda^+_\beta,\theta^+}$ we have $M_k \models \varphi[\bar
a^1] \Leftrightarrow M_\omega \models \varphi[\bar a^2]$.
\end{enumerate}
\end{enumerate}
\mn
The case $\bold s = \bold s^0_\beta$, the initial state, suffice to prove
the desired result.
\medskip

\noindent
\underline{First Case}:  $\beta = 0$.

Trivial.
\medskip

\noindent
\underline{Second Case}:  $\beta$ a limit ordinal.

If the AIS makes its move choosing $(A,\iota,\beta_1)$, the ISO player may
pretend $\beta_{\bold s} = \beta +1$ and use the induction
hypothesis, this is O.K. as $\beta < \beta +1$, as in proving \ref{a14}.
\medskip

\noindent
\underline{Third Case}:  $\beta = \gamma +1$.

As in the proof of \ref{a14}.
\end{PROOF}

We can sum up the easy properties but first we present two
definitions.
\begin{definition}
\label{a21k}
1) We say that $N$ code $[N]^{\le \lambda}$ \when \, for some two-place
predicate $R$ we have $\{\text{set}(b,N):b \in N\}$ is a cofinal
subset of $([N]^{\le \lambda},\subseteq)$.

\noindent
2) In this case for $b \in N$ let set$(b,N) = \text{ set}_R(b,N) =
\{a:a R^N b\}$ and if $R$ is clear from the context we may omit it. 
\end{definition}

\noindent
For our results it seems helpful to define some variants of ``$\cL$
satisfies the downward LST".
\begin{definition}
\label{a19}
Let $\cL$ be a logic, $\tau_0$ a vocabulary $\psi_0 \in \cL(\tau_0)$
a sentence with $\tau_0$ of minimal cardinality \underline{or} $T \subseteq
\cL(\tau_0)$, i.e. a set of sentences with $|\tau_0|$ minimal.

\noindent
1) $\cL$ satisfies LST$^\iota_{< \kappa}$ \when \, for every
   vocabulary $\tau$ of cardinality $< \kappa$ and $\varphi \in
   \cL(\tau)$ we have LST$^\iota_{< \kappa}(\varphi)$, see below;
   if $\iota =0$ we may omit it.

If $\kappa = \lambda^+$ we may write ``$\lambda$" instead ``$<
\kappa$".  LST$^\iota_{< \kappa}(\psi)$ means LST$^\iota_{<
\kappa}(\{\psi\})$ similarly below and LST$^\iota_{< \kappa}(\le
\theta)$ means LST$^\iota_{<\kappa}(T)$ when $|T| \le \theta$; the
reader may concentrate on the case $\kappa = \lambda^+$.

\noindent
2) For $\iota=0$, let LST$^\iota_{<\kappa}(T)$ means: if $T$ has a
   model $N$ of cardinality $\ge \kappa$ \then \, it has a model $M$
   of cardinality $\lambda$ for arbitrary large cardinals $\lambda < \kappa$.

\noindent
3) If $\iota = 1$, similarly but $M \subseteq N$ (moreover can demand $A
 \subseteq M$ for a given $A \in [N]^{<\kappa}$.

\noindent
4) For $\iota=2$, let LST$^\iota_{<\kappa}(T)$ means: if $N$ is a
   model of $T$ of cardinality $\ge \kappa$ \then \, for arbitrarily
   large $\lambda < \kappa$ we have: if $N^+$ is an expansion of $N$
   and $\tau(N^+)$ has cardinality $\le \lambda$ \then \, some $M^+
\subseteq N^+$ of cardinality $\lambda$ satisfies $M = M^+ \rest \tau_T$
 is a model of $T$.

\noindent
5) For $\iota=3$, as in $\iota=2$, but we further assume that $N^+$
   code $[N]^{\le \lambda}$, see \ref{a21k} below and we further conclude that
   for some $\langle d_n:n < \omega\rangle \in {}^\omega(M^+)$ we have
   $|M^+| = \cup\{\text{set}(d_n,N^+):n < \omega\}$.

\noindent
6) We define LST$^{\iota,*}_{< \kappa}(T)$ similarly but for a
   pregiven $\mu < \kappa$ can demand $\|M\| = \|M\|^\mu$ and similarly
   LST$^{\iota,*}_{< \kappa}$.
\end{definition}

\begin{conclusion}
\label{a21}
Assume $\theta < \kappa = \beth_\kappa$.

\noindent
1) $\bbL^1_{\le \theta}$ is a nice logic with full substitution, see
Definition \ref{z29}.

\noindent
2) $\bbL^1_\kappa$ is a nice logic and if $\kappa$ is regular
   (equivalently strongly inaccessible) then it has full subsitution.

\noindent
3) $\bbL^{-1}_\kappa \le \bbL^1_\kappa \le \bbL^0_\kappa$, in
   fact $\bbL_{\theta^+,\aleph_0} \le \bbL^1_{\le \theta} \le
\bbL_{\beth_{\theta^+},\theta^+}$ (really $\bbL^1_{\le \theta} \le
   \bbL_{\beth_{\theta^+},\theta^+,\theta^+}$ and even $\le
   \bbL_{\infty,\theta^+,\theta^+}$).   

\noindent
4) $\bbL^1_\kappa$ satisfies the following versions of the downward LST:
\mn
\begin{enumerate}
\item[$(a)$]   every sentence $\psi \in \bbL^1_\kappa$ which has a
model has a model of cardinality $< \kappa$
\sn
\item[$(b)$]  $\bbL^1_\kappa$ has LST$^\iota_{< \kappa}$ for $\iota
\le 3$ (see \ref{a19} above)
\sn
\item[$(c)$]  for every $\psi \in \bbL^1_\kappa$ there is $\partial <
\kappa$ such that: if $N$ is a model of $\psi$ of cardinality $\lambda$
and $\mu = \mu^{< \partial} \le \lambda$ or at least $\mu =
\sum\limits_{n < \omega} \mu_n \le \lambda$ and $(\mu_n)^{< \partial}
= \mu_n$ for $n < \omega$ \then \, some submodel $M$ of $N$ of
cardinality $\mu$ is a model of $\psi$.
\end{enumerate}
\end{conclusion}

\begin{PROOF}{\ref{a21}}
1) By \ref{a13m}.

\noindent
2) Follows as $\bbL^1_{\le \theta}$ is $\subseteq$-increasing with
  $\theta$ and the defintion of $\bbL^1_\kappa$ in Definition
\ref{a13}(2).

\noindent
3) We have $\bbL^1_{\le \theta} \le \bbL_{\beth_{\theta^+},\theta^+}$
 (and moreover $\bbL^1_{\le \theta} \le
   \bbL_{\beth_{\theta^+},\theta^+,\theta^+}$) by \ref{a14}.  So by the
   definition, $\bbL^1_\kappa \le \bbL^0_\kappa$.
Note that if $M_1,M_2$ are
   $\cE^1_{\text{qf}(\tau),\aleph_0,\alpha}$-equivalent then $M_1
   \models \psi \Leftrightarrow M_2 \models \psi$ for every $\psi \in
   \bbL_{\beth^+_\alpha,\aleph_0,\alpha}$ hence
   $\bbL_{\theta^+,\aleph_0} \le \bbL^1_{\le \theta}$ and so by the
   definitions, $\bbL^{-1}_\kappa \le \bbL^1_\kappa$.

\noindent 
4) Clause (a) follows by part (3) and the parallel results on
   $\bbL_{\lambda^+,\theta^+}$ of Hanf (see, e.g. \cite{Dic85}).  
Clause (b) follows by \ref{a18}.
\end{PROOF}
\newpage

\section {Serious properties of $\bbL^1_\kappa$}

First we prove a strong form of non-definability of well ordering.
\begin{claim}
\label{a22}
Let $\kappa = \beth_\kappa$ so a strong limit cardinal.

\noindent
1) Property (d) of Problem \ref{z2} holds.

\noindent
2) Moreover, if (A) then (B) \underline{where} $\tau$ is a vocabulary
to which the predicates $P,<,R$ (unary, binary, binary) belongs and:
\mn
\begin{enumerate}
\item[$(A)$]  $(a) \quad \psi \in \bbL^1_{\le \theta}(\tau)$ is
defined using $(\text{\rm qf}(\tau),\theta,\alpha)$
\sn
\item[${{}}$]  $(b) \quad \theta \ge |\tau|$
\sn
\item[${{}}$]  $(c) \quad \partial = \beth_{\alpha +1}(\theta)$
\sn
\item[${{}}$]  $(d) \quad \mu = 2^\partial$
\sn
\item[${{}}$]  $(e) \quad M \models \psi$
\sn
\item[${{}}$]  $(f) \quad (P^M,<^M) \cong ((2^\mu)^+,<)$
\sn
\item[${{}}$]  $(g) \quad M$ is a $\tau$-model satisfying:

\hskip25pt $\{\{a:a R^M b\}:b \in M\} = [M]^{\le \mu}$ 

\hskip25pt or just the former family is cofinal in the latter family
(both ordered 

\hskip25pt by inclusion)
\sn
\item[$(B)$]  there are a $\tau$-model $N$ and a sequence
$\langle b_n:n < \omega\rangle$ such that:
\begin{enumerate}
\item[$(\alpha)$]  $N \models \psi$
\sn
\item[$(\beta)$]  $(P^N,<^N)$ is not well ordered
\sn
\item[$(\gamma)$]  $b_n \in N$ and $N \models (\forall x)(x
R b_n \rightarrow x R b_{n+1})$
\sn
\item[$(\delta)$]  $N \models (\forall b) \bigvee\limits_{n < \omega} 
[b R b_n]$.
\end{enumerate}
\end{enumerate}
\end{claim}

\begin{PROOF}{\ref{a22}}
1) Follows by (2).

\noindent
2) Without loss of generality 
$|M|$ is an ordinal, $P^M = (2^\mu)^+$ and $<^M$ is the 
usual order of the ordinals on $|M|$.
By induction on $n$ we choose a sequence $\langle
(M_{n,\gamma},\bar\beta_{n,\gamma},\bar b_{n,\gamma}):\gamma <
(2^\mu)^+\rangle$ such that
\mn
\begin{enumerate}
\item[$\boxplus_1$]  $(a) \quad \bar\beta_{n,\gamma} = \langle
\beta_{n,\gamma,\ell}:\ell < n\rangle$ is a decreasing sequence of ordinals
from $P^M$,

\hskip25pt  which are $> \gamma$ and $\bar b_{n,\gamma} = \langle
b_{n,\gamma,\ell}:\ell < n \rangle$ is a sequence 

\hskip25pt of members of $M$
\sn
\item[${{}}$]  $(b) \quad M_{n,\gamma} \prec_{\bbL_{\partial^+,\theta^+}}M,
\|M_{n,\gamma}\| = \mu$ and $b_{n,\gamma,\ell},\beta_{n,\gamma,\ell} 
\in M_{n,\gamma}$ for $\ell < n$

\hskip25pt (recall $\mu = 2^\partial$)
\sn 
\item[${{}}$]  $(c) \quad$ for $\gamma_1 < \gamma_2 < (2^\mu)^+$ we have
$(M_{n,\gamma_1},\bar b_{n,\gamma_1},\bar\beta_{n,\gamma_1}) \cong
(M_{n,\gamma_2},\bar b_{n,\gamma_2},\bar\beta_{n,\gamma_2})$, 

\hskip25pt  note that the isomorphism 
is unique as $<^M$ is a well ordering and 

\hskip25pt necessarily the
isomorphism maps $b_{n,\gamma_1,\ell}$ to $b_{n,\gamma_2,\ell}$ 

\hskip25pt and $\beta_{n,\gamma_1,\ell}$ to $\beta_{n,\gamma_2,\ell}$ 
for $\ell < n$
\sn
\item[${{}}$]  $(d) \quad$ if $n = m+1$ then
\begin{enumerate}
\item[${{}}$]  $(i) \quad \bar\beta_{n,\gamma} \rest m =
\bar\beta_{m,\zeta}$ when $\zeta = \beta_{n,\gamma,m}$
\sn
\item[${{}}$]  $(ii) \quad M_{m,\beta_{n,\gamma,m}} 
\prec_{\bbL_{\partial^+,\theta^+}} M_{n,\gamma}$
\sn
\item[${{}}$]  $(iii) \quad a \in 
M_{m,\beta_{n,\gamma,m}} \Rightarrow M \models a R b_{n,\gamma}$
\sn
\item[${{}}$]  $(iv) \quad$ if $b \in M_{m,\beta_{n,\gamma,m}}$ then
$\{a:M \models a Rb\} \subseteq \{a:M \models a R b_{n,\gamma}\}$.
\end{enumerate}
\end{enumerate}
\mn
For $n=0$ just choose $M_{0,\gamma} = M_0 \prec_{\bbL_{\partial^+,\theta^+}}
M,\|M_0\| = \mu$, (OK as $\mu^\partial = \mu$, in fact an overkill).

For $n=m+1$, choose $M'_{m,\gamma}, b'_{m,\gamma}$ such that:
\mn
\begin{enumerate}
\item[$\boxplus_2$]  $(i) \quad$ if $a \in M_{m,\gamma}$ then $a
R^M b'_{m,\gamma}$ 
\sn
\item[${{}}$]  $(ii) \quad$ moreover if $b \in M_{m,\gamma}$ then $\{a
\in M:M \models a R b\} \subseteq$

\hskip25pt $\{a \in M:M \models a R b'_{m,\gamma}\}$
\sn
\item[${{}}$]  $(iii) \quad M_{m,\gamma} 
\cup \{b'_{m,\gamma}\} \cup \{\gamma\} \subseteq M'_{m,\gamma}$
\sn
\item[${{}}$]  $(iv) \quad M'_{m,\gamma} \prec_{\bbL_{\partial^+,\theta^+}} M$.
\end{enumerate}
\mn
The number of $(M'_{m,\gamma},M_{m,\gamma},\bar\beta_{m,\gamma
+1},\bar b_{n,\gamma +1},\gamma,b'_{m,\gamma})/\cong$ is $\le
2^\mu$ so for some unbounded $Y_n \subseteq (2^\mu)^+$, the models
$\langle (M'_{m,\gamma},M_{m,\gamma +1},\bar\beta_{m,\gamma +1},
\bar b_{n,\gamma +1},\gamma,b'_{m,\gamma}):\gamma \in Y_n\rangle$ are
pairwise isomorphic.

For any $\zeta < (2^\mu)^+$ let $\gamma_\zeta = \text{ Min}(Y
\backslash (\zeta +1)),M_{n,\zeta} = M'_{m,\gamma_\zeta},\bar b_{\zeta,n} =
\bar b_{m,\gamma_\zeta} \char 94 \langle b'_{m,\gamma_\zeta}\rangle,
\bar\beta_{n,\zeta} = \bar\beta_{m,\gamma_\zeta} 
\char 94 \langle \gamma_\zeta \rangle$.  Clearly $\langle
(M_{n,\zeta},\bar b_{n,\zeta},\bar\beta_{n,\zeta}):\zeta <
(2^\mu)^+\rangle$ is as required.

Having carried the induction, it is easy to find models 
$N_n$ and elements $a_\ell,b_\ell$ for $\ell < n$ by induction on $n <
\omega$ such that:
\mn
\begin{enumerate}
\item[$\boxplus_3$]  $(a) \quad
(N_n,a_0,\dotsc,a_{n-1},b_0,\dotsc,b_{n-1})$ is isomorphic to 
$(M_{\gamma,n},\bar \beta_{\gamma,n},\bar b_{\gamma,n})$ 

\hskip25pt for every $\gamma < (2^\mu)^+$
\sn
\item[${{}}$]  $(b) \quad N_m \prec_{\bbL_{\partial^+,\theta^+}} N_n$ if $m<n$.
\end{enumerate}
\mn
Now easily by \ref{a18} the model $N = 
\bigcup\limits_{n} N_n$ and the sequence $\bar b = 
\langle b_n:n < \omega\rangle$ are as required in clause (B) 
with $\langle a_n:n < \omega\rangle$ witnessing clause $(\beta)$ of
clause $(B)$.
\end{PROOF}

\begin{conclusion}
\label{a24}
For $\kappa = \beth_\kappa$, the logic $\bbL^1_\kappa$ satisfies:

\noindent
\underline{SUDWO}$^1_\kappa$ \,\, (strong undefinability of well
ordering) which means:

if $\psi \in \bbL^1_\kappa(\tau),|\tau| < \kappa$ and $<,R$ 
are two place predicates from $\tau$ \then \,
for every large enough $\mu_1 < \kappa$ for arbitrarily large enough 
$\mu_2 < \kappa$ we have:
\mn
\begin{enumerate}
\item[$\circledast$]  \If \, $\lambda > \mu_2$ and $\gA$ 
is a $\tau$-expansion of $(\cH(\lambda),\in,\mu_1,\mu_2,<)$ 
with $<$ the order on the ordinals and $R^{\gA}$ being $\in$, that
is $\in \rest \cH(\lambda)$ \then \, we can find $\gB,a_n,d_n$
(for $n < \omega$) such that
\sn
\item[${{}}$]  $(a) \quad \gB \models \psi \Leftrightarrow M \models
\psi$
\sn
\item[${{}}$]  $(b) \quad \gB \models ``d_{n+1} < d_n < \mu_2"$ for $n
< \omega$
\sn
\item[${{}}$]  $(c) \quad \gB \models ``a_n \subseteq a_{n+1}$ has
cardinality $\le \mu_1"$
\sn
\item[${{}}$]  $(d) \quad$ if $e \in \gB$ then $\gB \models 
``e \in a_n"$ for some $n$

\qquad we may add
\sn
\item[${{}}$]  $(d)^+ \quad$ if $e \in \gB$ and $\gB \models ``|e| \le
\mu_1"$ \then \, $\gB \models ``e \subseteq a_n"$ for some $n$.
\end{enumerate}
\end{conclusion}

\begin{remark}
\label{a26}
There are other variants.  At the moment the distinction is not
crucial.
\end{remark}

\begin{PROOF}{\ref{a24}}
By \ref{a22}.
\end{PROOF}

\begin{theorem}
\label{a28}
\underline{First Characterization Theorem}

We have $\cL \le \bbL^1_\kappa$, i.e. 
every $\psi \in \cL$ is equivalent to some $\psi' \in \bbL^1_\kappa$
(i.e. for any vocabulary $\tau,\ldots$) \when \, :
\mn
\begin{enumerate}
\item[$\boxplus$]  $(a) \quad \kappa = \beth_\kappa$ hence
is a strong limit uncountable cardinal 
\sn
\item[${{}}$]  $(b) \quad \cL$ is a nice logic 
\sn
\item[${{}}$]  $(c) \quad \Upsilon$, the 
occurance number of $\cL$, is $\le \kappa$
\sn
\item[${{}}$]  $(d) \quad 
\bbL^{-1}_\kappa \le \cL$, i.e. $\bbL_{\theta^+,\aleph_0}
\le \cL$ for $\theta < \kappa$ or just
\sn
\item[${{}}$]  $(d)^- \quad$ for every ordinal $\alpha < \kappa$, the model
$M_\alpha$ expanding $(\alpha,<)$ is

\hskip25pt  characterized up to
isomorphism by some sentence from $\bbL^1_\kappa(\tau(M_\alpha))$, 

\hskip25pt so $|\tau(M_\alpha)| < \kappa$
\sn
\item[${{}}$]  $(e) \quad \bbL$ satisfies {\rm SUDWO}$^1_\kappa$.
\end{enumerate}
\end{theorem}

\begin{remark}
\label{a29}
1) We could add the downward LST and restrict somewhat SUDWO$_\kappa$.

\noindent
2) As $\bbL^1_\kappa$ satisfies the demands (by \ref{a21}) we can
   rephrase Theorem \ref{a28} as: $\bbL$ is $\le$-maximal such that
   (a)-(e) holds.
\end{remark}

\begin{PROOF}{\ref{a28}}
Assume toward contradiction that $\psi^* \in
\cL(\tau)$ is a sentence which is not equivalent to any sentence of
$\bbL^1_\kappa(\tau)$.  As the occurance number of $\cL$ is $\le
\kappa$ by clause (c) of the assumption, \wilog \, the vocabulary
$\tau$ has cardinality $< \Upsilon \le \kappa$.  Similarly by $\cL$
being a logic, using clauses (c),(d) of Definition \ref{z21}, without
loss of generality the symbols we add to $\tau$ in $(*)_3(a)$ below do
not belong to $\tau$.
We shall derive another sentence $\psi^{**}$ in a somewhat bigger
vocabulary, and let $(\mu_1,\mu_2)$ be as in Definition \ref{a24} for
$\psi^{**}$ and $\mu_1 \ge \theta_*$.

Let
\mn
\begin{enumerate}
\item[$(*)_0$]  $(a) \quad \{P_i:i < i(*)\}$ list the predicates in
$\tau$
\sn
\item[${{}}$]  $(b) \quad \{H_j:j < j(*)\}$ list the function symbols
in $\tau$
\sn
\item[${{}}$]  $(c) \quad$ let the vocabulary $\tau_1$ be
\sn
\begin{enumerate}
\item[${{}}$]  $\bullet \quad$ the set of predicates $\{P'_i:i \le i(*)\}$
\sn
\item[${{}}$]  $\bullet \quad$ the set of function symbols 
$\{H'_j:j < j(*)\}$
\sn
\item[${{}}$]  $\bullet \quad$ arity$_{\tau_1}(P'_{i(*)})=2$
\sn
\item[${{}}$]  $\bullet \quad$ arity$_{\tau_1}(P'_i) = \text{
arity}_\tau(P_i)+1$ for $i<i(*)$
\sn
\item[${{}}$]  $\bullet \quad$ arity$_{\tau_1}(H'_j) = \text{
arity}_\tau(H_j)+1$ for $j<j(*)$
\end{enumerate}
\item[${{}}$]  $(d) \quad$ for a $\tau_1$-structure $\gB$ and $c \in
\gB$ let $N^{\gB}_c$ be the $\tau$-model (if exists) 

\hskip35pt such that
\sn
\begin{enumerate}
\item[${{}}$]  $\bullet \quad$ it has universe
$\{d:\gB \models P'_{i(*)}(d,d)\}$
\sn
\item[${{}}$]  $\bullet \quad (P_i)^{N^{\gB}_c} =
\{\bar a:\langle c \rangle \char 94 \bar a \in (P'_i)^{\gB}\}$
for $i < i(*),j <j(*)$
\sn
\item[${{}}$]  $\bullet \quad H^{N^{\gB}_c}_j(\bar a) =
(H'_j)^{\gB}(\langle c \rangle \char 94 \bar a)$ for $i<i(*)$.
\end{enumerate}
\end{enumerate}
\mn
For each $\alpha < \kappa$, let $\theta(\alpha) = \theta_\alpha 
= |\alpha| + \aleph_0$, so the
equivalence relation $\cE^1_{\text{qf}(\tau),\theta(\alpha),\alpha}$ 
is well defined, and has $< \beth_{\theta(\alpha)^+}$ equivalence classes 
say $\langle N_{\alpha,\varepsilon}/\cE^1_{\text{qf}(\tau),
\theta(\alpha),\alpha}:\varepsilon < \varepsilon_\alpha 
< \beth_{\theta(\alpha)^+}(|\tau|)\rangle$. 

By the definition of $\bbL^1_{\le  \theta}(\tau)$ which is $\le
\bbL^1_\kappa(\tau)$, for each such 
pair $(\alpha,\varepsilon)$ there is a sentence
$\vartheta_{\alpha,\varepsilon} \in \bbL^1_{\le \theta(\alpha)}(\tau)$
which define
$N_{\alpha,\varepsilon}/\cE^1_{\text{qf}(\tau),\theta(\alpha),\alpha}$
 and moreover for every $u \subseteq \varepsilon_\alpha$
the sentence $\vartheta_{\alpha,u} =
\vee\{\vartheta_{\alpha,\varepsilon}:\varepsilon \in u\}$ belongs to
$\bbL^1_\kappa(\tau)$, i.e. up to equivalence.  
Hence by our assumption toward contradiction for
some $\zeta(\alpha) < \varepsilon_\alpha$ there are
\mn
\begin{enumerate}
\item[$(*)_1$]  $(a) \quad M_\alpha,N_\alpha \in
N_{\alpha,\zeta(\alpha)}/\cE^1_{\text{qf}(\tau),\theta(\alpha),\alpha}$
\sn
\item[${{}}$]  $(b) \quad N_\alpha \models \psi^*$ but $M_\alpha
\models \neg \psi^*$.
\end{enumerate}
\mn
By the definition of $\cE^1_{\text{qf}(\tau),\theta(\alpha),\theta}$
there is a sequence $\bar M_\alpha$ such that
\mn
\begin{enumerate}
\item[$(*)_2$]  $(a) \quad \bar M_\alpha = \langle M_{\alpha,k}:k \le
k(\alpha)\rangle$ so $k(\alpha) < \omega$
\sn
\item[${{}}$]  $(b) \quad M_{\alpha,0} = M,M_{\alpha,k(\alpha)} =
N_\alpha$
\sn
\item[${{}}$]  $(c) \quad M_{\alpha,k},M_{\alpha,k +1}$ are
$\cE^0_{\text{qf}(\tau),\theta(\alpha),\alpha}$-equivalent for
$k<k(\alpha)$
\sn
\item[${{}}$]  $(d) \quad$ let $M_{\alpha,k} \cong M_{\alpha,k(\alpha)}$
for $k \in (k(\alpha),\omega)$.
\end{enumerate}
\mn
Without loss of generality the universes of the models
$M_{\alpha,k} (\alpha < \kappa,k < \omega)$
are pairwise disjoint and $\tau \subseteq \cH_{< \aleph_0} (|\tau|)$.
For large enough $\lambda$ we have $\langle \langle M_{\alpha,k}:k \le
k(\alpha)\rangle:\alpha < \kappa\rangle \in \cH(\lambda)$ for every
$\alpha < \kappa$ and so $\lambda > \kappa$; 
let $\gA_\lambda = \gA(\lambda)$ be a model such that:
\mn
\begin{enumerate}
\item[$(*)_4$]  $(a) \quad$ the vocabulary $\tau_2$ of $\gA_\lambda$ is 
$\tau_1 \cup \{<,E,R,R_1,F_1,F_2\} \cup \{c_\alpha:\alpha \le 
\theta^+_*\}$,

\hskip25pt  with $c_\alpha$ individual constants,
$F_1,F_2$ unary function symbols, 

\hskip25pt $<,E$ binary predicates
\sn
\item[${{}}$]  $(b) \quad$ universe $\cH(\lambda)$
\sn
\item[${{}}$]  $(c) \quad$ among the elements $c^{\gA(\lambda)}_\alpha,
(\alpha < \theta_*)$ will be $\kappa,\theta_*,\theta^+_*,\mu_1,
\mu_2,\tau,\tau_1$ and 

\hskip25pt  every symbol in $\tau$ and $\alpha \le \theta_* 
\Rightarrow c_{\theta^* +\alpha}(\gA) = \alpha$; so

\hskip25pt $c^{\gA(\lambda)}_0 = \kappa,c^{\gA(\lambda)}_1 =
\theta_*$, etc.
\sn
\item[${{}}$]  $(d) \quad$ the functions $\omega \alpha +k \mapsto
M_{\alpha,k}$ and $\alpha \mapsto k(\alpha)$ are 
$F^{\gA(\lambda)}_1,F_2^{\gA(\lambda)}$ 

\hskip25pt respectively
\sn
\item[${{}}$]  $(e) \quad <^{\gA(\lambda)}$ is the order on the ordinals
\sn
\item[${{}}$]  $(f) \quad E^{\gA(\lambda)}$ is an equivalence 
relation such that

\hskip25pt the equivalence classes of $E^\gA$ are 
$\{|M_{\alpha,k}|:\alpha < \kappa$ and $k < \omega\}$
\sn
\item[${{}}$]  $(g)(\alpha) \quad R^{\gA(\lambda)} = \{(a,b):a \in b \in
\cH(\lambda)$ and $|b| \le \mu_1\}$
\sn
\item[${{}}$]  $\quad(\beta) \quad R^{\gA(\lambda)}_1 = \theta \rest 
\cH(\lambda)$
\sn
\item[${{}}$]  $(h) \quad (P'_{i(*)})^{\gA(\lambda)} = 
\{(a,\omega \alpha + k):k < \omega,\alpha < \kappa$ and $a \in
M_{\alpha,k}\}$
\sn
\item[${{}}$]  $(i) \quad (P'_i)^{\gA(\lambda)} = \{\bar a \char 94 
\langle \omega \alpha + k \rangle:\alpha < \kappa,\bar a \in
P^{M_{\alpha,k}}_i$ and $k < \omega\}$ for $i < i(*)$
\sn
\item[${{}}$]  $(j) \quad H^{\gA(\lambda)}_j$ is an
$(\text{arity}_\tau(H_j) + 1)$-place function satisfying

\hskip25pt $H^{\gA}_j(\bar a \char 94 \langle \omega \alpha + k \rangle) =
H^{M_{\alpha,k}}_j(\bar a)$ for $j < j(*),\alpha < \kappa$ and $k < \omega$.
\end{enumerate}
\mn
Let $\psi^{**}$ says:
\mn
\begin{enumerate}
\item[$(*)_5$]  $(a)\quad$ one first order sentence saying all 
relevant (set-theoretic) properties
\sn
\item[${{}}$]  $(b) \quad$ description of $\tau,\tau_1$ 
\sn
\item[${{}}$]  $(c) \quad (\forall \alpha < \kappa)(M_{\alpha,0} \models
\psi^*$ and $M_{\alpha,k(\alpha)} \models \neg \psi^*)$.
\end{enumerate}
\mn
Why possible?  For clauses (a),(b) as $\bbL^{-1}_\Upsilon
\le \cL$, by clause (d) of the assumption or just ensuring $\{a:a <
c_1\}$ isomorphic to $\theta^+_*$ if we just assume clause (d)$^-$ of the
assumption and then the rest is said by a first order sentence.  
For clause (c) as the logic has restriction and weak
substitution (see clauses (c),(d) of Definition \ref{z29}) recalling the logic 
is closed under the first order operations.

Now apply SUDWO$^1_\kappa$ to $\psi^{**},(\gA_\lambda,\mu_1,\mu_2)$,
and we get $\gB,a_n,d_n (n < \omega)$ as there.

Now
\mn
\begin{enumerate}
\item[$(*)_6$]  \wilog \, $\tau,\tau_1$ are interpreted in $\gB$ as
$\tau,\tau_1$ respectively;
 similarly the individual constant $c_\alpha$ for $\alpha <
\theta^+_*$ is interpreted as $c^{\gA(\lambda)}_\alpha$,
hence if $\gB \models$ ``$c$ is an ordinal $< \kappa$ and 
$k(c) = F_2(c)$" \then \, for $k < k(c),m \le k(c)$ we have:
\begin{enumerate}
\item[$(a)$]  $M^{\gB}_{c,m}$ defined naturally as $N^{\gB}_{F^{\gB}_1(c,m)}$
 is a $\tau$-model (i.e. no non-standard predicates!)
\sn
\item[$(b)$]   $M^{\gB}_{c,0} \models \psi$ and $M^{\gB}_{c,k(c)} 
\models ``\neg \psi"$ 
\sn
\item[$(c)$]  $\gB \models$ ``in the game $\Game_{\text{qf}(\tau),\theta_*,c}
[M^{\gB}_{c,k},M^{\gB}_{c,k+1}]$ the ISO player wins"
\end{enumerate}
\item[$(*)_7$]   let $k(*) = F^{\gB}_2(d_0)$
\sn
\item[$(*)_8$]  we fix a winning strategy {\bf st}$_k$ for ISO in 
$\gB$'s sense in the game
$\Game_{\qf(\tau),\theta_*,d}[N^{\gB}_{d_0,k},N^{\gB}_{d_0,k+1}]$.
\end{enumerate}
\mn
Now for each $k < k(*)$ by induction on $n < \omega$ we 
choose $\bold s^{n,k} \in \gB$ such that:
\mn
\begin{enumerate}
\item[$(*)_9$]  $(a) \quad \gB \models ``\bar{\bold s}^{n,k} = \langle \bold
s_{0,k},\dotsc,\bold s_{n,k}\rangle$ is an initial segment of a play of the
game 

\hskip25pt $\Game_{\qf(\tau),\theta_*,d_0}
[M^{\gB}_{d_0,k},M^{\gB}_{d_0,k+1}]"$
\sn
\item[${{}}$]  $(b) \quad \gB \models$ ``in this initial segment
$\bold s^{n,k}$ the ISO player uses his 

\hskip25pt winning strategy {\bf st}$_k$"
\sn
\item[${{}}$]  $(c) \quad$ the AIS player chooses:
\begin{enumerate}
\item[${{}}$]  $(\alpha) \quad$ if $n$ is even \then \, $\iota_n = 1$ and 
$A'_n = \{e \in \gB:e \in M^{\gB}_{d_0,k}$

\hskip35pt  that is $\gB \models ``P'_{i(*)}(e,F_1(d_0,k))"$ and
$\gB \models ``e R a_n"$
\sn
\item[${{}}$]  $(\beta) \quad$ if $n$ is odd \then \, $\iota_n = 2$ and
$A'_n = \{e \in \gB:e \in M^{\gB}_{d_0,k+1}$ and

\hskip35pt  $\gB \models ``e R a_n"\}$
\sn
\item[${{}}$]  $(\gamma) \quad$ the ``ordinal" $\beta_{\bold s_{n,k}}=d_{n+1}$.
\end{enumerate}
\end{enumerate}
\mn
This can be done and $g_k = \{(a_1,a_2)$: for some $n$ we have
$\gB \models ``g^{\bold s_{n,k}}_n(a_1) = a_2"\}$, it is 
an isomorphism from $M^{\gB}_{d_0,k}$ onto $M^{\gB}_{d_0,k+1}$.  As
this holds for every $k < k(*)$ we get that
$M^{\gB}_{d_0,0},M^{\gB}_{d_0,k(*)}$ are isomorphic so by $(*)_6(b)$ we
get a contradiction.
\end{PROOF}

\begin{definition}
\label{a24}
Let $\Theta_{\theta,R}$ be the sentence (for $R$ a binary predicate)
such that $M = (|M|,R^M) \models \Theta_{\theta,R}$ iff
$\cP = \{\{a:aR^Mb\}:b \in M\}$ is a $(\theta,\aleph_0)$-cover of
$\Dom(P^M)$ which means: it is a family of subsets of $\Dom(R) =
\{a:a R^M b$ for some $b \in M\}$ each of cardinalilty
$\le \theta$ such that any $u \in [\Dom(R^M)]^{\le \theta}$ is included in the
union of countably many such sets.
\end{definition}

\begin{claim}
\label{a30}
1) $\Theta_{\theta,R}$ can be expressed by a sentence in 
$\bbL^1_{\le \theta^+}(\{R\})$.

\noindent
2) In fact, if $M_1,M_2$ are $\cE^1_{\text{\rm qf}(\{R\}),\theta^+,\omega
   + \omega +1}$-equivalent $\{R\}$-models \then \, $M_1 \models
   \Theta_{\theta,R} \Leftrightarrow M_2 \models \Theta_{\theta,R}$.

\noindent
3) If $M_1 \cE^1_{\qf(\tau),\theta,\alpha} M_2$ and $\alpha \ge
   \omega + \omega +1$ and {\rm cf}$(\theta) \ge \aleph_0$ \then \,
   $\|M_1\| = \|M_1\|$ or $\|M_1\|,\|M_2\| > \theta$.

\noindent
4) If $M_1 \cE^1_{\qf(\tau),\theta,\alpha} M_2$ and $\|M_1\| \le
   \theta$ and $\alpha \ge \omega + \omega +1$ \then \, $M_1 \cong M_2$.
\end{claim}

\begin{PROOF}{\ref{a30}}
1) By (2).

\noindent
2) Toward contradiction, assume this fails; and \wilog \, $M_1
\cE^0_{\qf(\{R\}),\theta^+,\omega + \omega +1} M_2$.
By symmetry, \wilog \, assume $M_2 \models
   \Theta_{\theta,R}$ but $M_1 \models \neg \Theta_{\theta,R}$.  We
   simulate a play of the game 
$\Game_{\qf(\{R\}),\theta^+,\omega + \omega +1}(M_1,M_2)$ in which the
   ISO player uses a (fixed) winning strategy $\text{\bf st}$.
\medskip

\noindent
\underline{Case 1}:  There is $B_* \in [M_1]^{\le \theta}$ not included
in a countable union of sets from $\{\{a:aR^{M_1}b\}:b \in M_1\}$.
Here $\Game_{\text{qf}(\{R\}),\theta,\omega + \omega +1}[M_1,M_2]$
suffice; we now simulate a play in which the ISO player uses {\bf st}.

In the first move we let the AIS player choose $\beta_{\bold s_1} = \omega +
\omega$ and choose $\iota_{\bold s_1} = 1,
A^1_{\bold s_1} = B_* \in [M_1]^{\le \theta}$, see
above.  Let the ISO player (using {\bf st}) 
complete the choice of the state $\bold
s_1$ (which is a winning state for itself, of course).

Let $n(*)$ be minimal such that:
\mn
\begin{enumerate}
\item[$(*)$]  $B := \{a \in A^1_{\bold s_1}:h^1_{\bold s_1}(a) <
n(*)-7\}$ is not included in a countable union of sets from
$\{\{a:aR^{M_1} b\}:b \in B\}$.
\end{enumerate}
\mn
By the case assumption, $n(*)$ is well defined.
In the following $n(*)-2$ moves, the AIS player takes care that
$\beta_{\bold s_2} = \omega + n(*),\dotsc,\beta_{\bold
s_{n(*)}} = \omega +1$ (and $\iota_{\bold s_2} = \ldots 
= \iota_{\bold s_{n(*)}} = 1$).

So $g_{\bold s_{n(*)}}$ is a function whose domain includes $\{a \in
A^1_{\bold s_1}:h^1_{\bold s_1}(a) < n(*)\}$ hence includes $B$.  Now
as $M_2 \models \Theta_{\theta,R}$ we can find $c_n \in M_2$ for $n <
\omega$ such that $(\forall b \in B) \, 
\bigvee\limits_{n} g_{\bold s_{n(*)}}(b) R^{M_2} c_n$ and 
AIS player takes care that $\beta_{\bold s_{n(*)+1}} = \omega$ and
$\{c_n:n < \omega\} \subseteq A^2_{\bold s_{n(*)+1}}$.

The rest should be clear.
\medskip

\noindent
\underline{Case 2}:  There is $b_* \in M_1$ such that 
$B_* = \{a \in M_1:a R^{M_1} b_*\}$ has cardinality $> \theta$.  

Here $\Game_{\text{qf}(\{R\}),\theta^+,\omega +1}$ suffice, and let $B
\subseteq B_*$ be of cardinality $\theta^+$.  We simulate a play in
which the AIS player takes care that $\beta_{\bold s_1} =
\omega,\iota_{\bold s_1} = 1,A^1_{\bold s_1} = \{b_*\} \cup B$.

We let $n_*$ be minimal such that $h^1_{\bold s_1}(b_*) < n_*-7$ and
$|\{b \in B:h^1_{\bold s_1}(b) < n_*-7\}| \ge \theta^+$.

The rest should be clear.

\noindent
3) Should be clear.
\end{PROOF}

\noindent
Recall
\begin{definition}
\label{a31}
1) A logic $\cL$ is $\Delta$-closed \when \, : for vocabularies
   $\tau_1 \subseteq \tau_2$ with $\tau_2 \backslash \tau_1$ finite,
   and sentences $\psi,\vartheta \in \cL(\tau_2)$ \If \, $K_0 = \{M
   \rest \tau_1:M \models \psi\},K_1 = \{M \rest \tau_1:M \models
   \vartheta\}$ are complementary classes of $\tau_1$-models 
\then \, some $\varphi \in \cL(\tau_1)$ define $K_1$.

\noindent
2) A logic $\cL$ is strongly $\Delta$-closed \when \, : for relational
vocabularies
$\tau_1 \subseteq \tau_2$ with $\tau_2 \backslash \tau_1$ finite,
unary predicate $P \in \tau_2 \backslash \tau_1$ and sentences 
$\psi,\vartheta \in \cL(\tau_2)$, \If \, $K_0 := \{(M
\rest \tau_1) \rest P^M:M \models \psi$ and $P^M \ne \emptyset\}$ and $K_1 := 
\{(M \rest \tau_1) \rest P^M:M \models \vartheta$ and $P^M \ne
\emptyset\}$ are complementary classes of $\tau_1$-models 
\then \, some $\varphi \in \cL(\tau_1)$ define $K_1$.

\noindent
3) A logic $\cL$ has dullness-elimination \when \,: if $\tau_1,\tau_2$
 are relational vocabularies, i.e.  
with predicates only, $\tau_2 = \tau_1 \cup
\{P\},P$ a unary predicate, $\psi_2 \in \cL(\tau_2)$ and for every
$M \models \psi_2$ we have $Q \in \tau_1 \Rightarrow Q^M = Q^M \rest
   P^M$ and $P^M \ne \emptyset$ \then \, for some $\psi_1 \in
   \cL(\tau_1)$ we have $\{M:M$ a $\tau_1$-model of $\psi_1\} = \{(M
 \rest \tau_1) \rest P^M:M$ a $\tau_2$-model of $\psi_2\}$.

\noindent
4) We say the logic $\cL$ is $\Delta$-closed in $\bold C$ for $\bold
C$ a class of cardinals \when \,: if in part (1) we assume just that
$K_0 \cap K^\tau_{\bold C},K_0 \cap K^\tau_{\bold C}$ are 
complimentary in $K^\tau_{\bold C} = \{M:M$ a $\tau$-model of 
cardinality from $\bold C\}$ \then \, for some $\psi \in \cL(\tau_1)$ for
every $\tau_1$-model $M$ of cardinality $\in \bold C$ we have $M
\models \psi \Leftrightarrow M \in K_0$.
\end{definition}

\begin{claim}
\label{a31d}
1) The logic $\bbL^1_{\le \theta}$ has dullness-elimination.

\noindent
2) If $\tau$ is empty, $M_1,M_2$ are $\tau$-models of cardinality $>
   \theta$ \then \, $M_1,M_2$ are
 $\cE^1_{\text{\rm qf}(\tau),\theta,\alpha}$-equivalent for $\alpha <
   \theta^+$.

\noindent
3) Also $\|M_1\|,\|M_2\| \ge \theta,\cf(\theta) = \aleph_0$ is O.K. in
   part (2).
\end{claim}

\begin{PROOF}{\ref{a31d}}
1) Easy (or use the sum theorem (see \ref{a40}) and part (2)).

\noindent
2),3)  Easy.
\end{PROOF}

\begin{theorem}
\label{a32}
\underline{Second Characterization Theorem}

Let $\kappa = \beth_\kappa$, then $\bbL^1_\kappa$ is the minimal logic
$\cL$ (up to equivalence) such that:
\mn
\begin{enumerate}
\item[$\oplus$]  $(a) \quad \cL$ is a nice logic
\sn
\item[${{}}$]  $(b) \quad \cL^{-1}_\kappa \le \cL$,
i.e. $\bbL_{\theta^+,\aleph_0} \le \cL$ for $\theta < \kappa$
\sn
\item[${{}}$]  $(c) \quad \Theta_{\theta,R}$ is equivalent to some sentence
in $\cL$
\sn
\item[${{}}$]  $(d) \quad \cL$ is $\Delta$-closed
\sn
\item[${{}}$]  $(e) \quad \cL$ has dullness-elimination.
\end{enumerate}
\end{theorem}

\begin{remark}
\label{a32d}
Putting together \ref{a28}, \ref{a32} we get full characterization.
\end{remark}

\begin{PROOF}{\ref{a32}}
First $\bbL^1_\kappa$ satisfies clause (a) by \ref{a21}(2), clause
(b) by \ref{a21}, clause (c) by \ref{a30}(1) and clause (d) by Theorem
\ref{a35} and Observation \ref{a33}(1) below, and lastly $\bbL^1_\kappa$ has
dullness-elimination by Claim \ref{a31d}(1) so together 
$\bbL^1_\kappa$ satisfies
$\oplus$ of \ref{a32}. Second, we shall assume $\cL$ satisfies $\oplus$ and
$\psi \in \bbL^1_\kappa(\tau)$ and we shall find $\psi' \in \cL(\tau)$
equivalent to it (it will be $\varphi_1$ below).  By ``$\cL$ satisfies
weak substitution" \wilog \, $\tau$ has predicates only.

We deal with and define naturally 
$\psi^{[P]} \in \bbL^1_\kappa(\tau),P$ a new unary predicate and let $\tau' =
\tau \cup \{P\}$ so $\psi^{[P]}$ says $(M \rest \tau) \rest P^M$
satisfies $\psi$. 

Let $\theta,\alpha < \kappa$ and let $\psi \in \bbL^1_\kappa(\tau)$ 
be defined using
$(\text{qf}(\tau),\theta,\alpha)$, so $|\tau| < \kappa$.  We let
$\partial = \beth_{\alpha +1}(\theta + |\tau|)$ and $\mu =
2^\partial$.  Let $\tau^+ = \tau' \cup \{R_1,R\}$, arity$(R_1) = 3$,
arity$(R) = 2$ and $R,R_1 \notin \tau'$.  
For $i=0,1$ let $K_\ell$ be the class of $\tau^+$-models $M$ such that:
\mn
\begin{enumerate}
\item[$\boxplus_\ell$]  $(a) \quad (|M|,R^M) \models
\Theta_{\mu,R}$ and $P^M \ne \emptyset$
\sn
\item[${{}}$]  $(b) \quad \{B_b:b \in M\}$ is $\subseteq$-directed
where $B_b = \{a \in M:a R^M b\}$
\sn
\item[${{}}$]  $(c) \quad$ for each $b \in M$ the set
$\{(a_1,a_2):(a_1,a_2,b) \in R^M_1\}$ is a well ordering

\hskip25pt  of $B_b$ of order type $\le \mu$
\sn
\item[${{}}$]  $(d) \quad$ if $b_1,b_2 \in M$ and $B_{b_1} \subseteq
B_{b_2} \subseteq P^M$ then 

\hskip25pt $(M \rest \tau) \rest B_{b_1}
\prec_{\bbL_{\partial^+,\partial^+}} (M \rest \tau) \rest B_{b_2}$
in particular $B_{b_1}$ is 

\hskip25pt closed under $F^M$ for any function symbol 

\hskip25pt $F \in \tau$
\sn
\item[${{}}$]  $(e) \quad$ if $b \in M$, then for some $c$ we have
$B_c = B_b \cap P^M
\subseteq B_c \subseteq P^M$
\sn
\item[${{}}$]  $(f) \quad$ if $b \in M$ and $B_b \subseteq P$
\then \, $(M \rest \tau) \rest
B_b \models \psi^{[P]}$ iff $\ell=1$
\sn
\item[${{}}$]  $(g) \quad$ if $Q \in \tau$ then $Q^M = Q^M \rest P^M$.
\end{enumerate}
\mn
Now
\mn
\begin{enumerate}
\item[$(*)_1$]   clauses (b),(c),(d),(e),(g) of $\boxplus_\ell$ can be
described by some $\varphi \in \bbL_{(2^\mu)^+,\aleph_0}(\tau_1)$.
\end{enumerate}
\mn
[Why?  Note clause (c), so obvious.]
\mn
\begin{enumerate}
\item[$(*)_2$]   clause (f) can be expressed by a sentence from
  $\cL(\tau^+)$.
\end{enumerate}
\mn
[Why?  Let $\langle N_\alpha:\alpha < 2^\mu \rangle$ list the
  $\tau$-models with universe $\subseteq \mu$.  For each $b \in M$
  such that $B_b \subseteq P^M$ let $g_b$ be the unique one-to-one
  function from $B_b$ onto some ordinal $\gamma_b \le \mu$ such that
\mn
\begin{enumerate}
\item[$\bullet$]  if $a_1,a_2 \in B_b$ then $(a_1,a_2,b) \in R^M_1
  \Leftrightarrow g_b(a) < g_b(a_2)$.
\end{enumerate}
\mn
Hence for such $b$ there is a unique $\alpha_b$ such that $g_b$ is an
isomorphism from $(M \rest \tau) \rest B_b$ onto $N_\alpha$.  For each
$\alpha$ there is a formula $\psi_\alpha(x) \in
\bbL_{\mu^+,\aleph_0}(\tau_M)$ such that for every $b \in M$ we have
$M \models \psi_\alpha[b]$ \Iff \, $b$ is as above and $\alpha_b =
\alpha$.

So for every $u \subseteq 2^\mu,\psi'_u(x) = \bigvee\limits_{\alpha
  \in u} \psi_\alpha(x)$ belongs to
  $\bbL_{(2^\mu)^+,\aleph_0}(\tau_M)$ hence is equivalent to some
  $\psi''_u(x) \in \cL(\tau_M)$.  Now $(*)_2$ should be clear.]
\mn
\begin{enumerate}
\item[$(*)_3$]   $K_\ell$ is definable by some
$\varphi_\ell \in \cL(\tau^+)$ for $\ell = 1,2$.
\end{enumerate}
\mn
[Why?  By clause (b) of the assumption, some $\varphi' \in
\cL(\tau^+)$ is equivalent to $\varphi$ where $\varphi^*$ is from
$(*)_1$.  By clause (c) of the
assumption some $\varphi'' \in \cL(\tau^+)$ is equivalent to
$\Theta_{\theta,R}$.  But by clause (a) of the assumption $\cL$ is a
nice logic, hence (see clauses $(b)_2,(b)_3,(c)$ of Definition
\ref{z29}) there is a sentence $\varphi'_\ell \in \cL(\tau^+)$
equivalent to $\varphi_\ell$ from $(*)_3$.
By clause $(b)_2$ of Definition \ref{z29},
there is $\varphi_\ell \in \cL(\tau^+)$ equivalent to 
$\varphi' \wedge \varphi'_\ell$ so we are done.]
\mn
\begin{enumerate}
\item[$(*)_4$]   for $\ell=1,2$ the classes
$\{(M_1 \rest \tau) \rest P^{M_1}:M_1 \in K_1\},
\{(M_2 \rest \tau)P^M_i,M_2 \in K_2\}$ are equal to $\{M:M$
a $\tau$-model of $\psi\},\{M:M$ a $\tau$-model of $\neg \psi\}$,
respectively.
\end{enumerate}
\mn
By $\cL$ being $\Delta$-closed and having dullness-elimination (see
clauses (d),(e) of the theorem's assumption) we are done.
\end{PROOF}

\noindent
We can note and recall
\begin{observation}
\label{a33}
1) If the logic $\cL$ satisfies interpolation \then \, it is $\Delta$-closed.

\noindent
2) $\cL$ is $\Delta$-closed and has dullness-elimination \Iff \,
   $\cL$ is strongly $\Delta$-closed.

\noindent
3) If $\cL$ satisfies interpolation with finitely many sorts (so the
interpolant mentions only the common sorts) \then \, $\cL$ is
   strongly $\Delta$-closed.

\noindent
4) The logic $\cL$ satisfies interpolation and dullness-elimination
   \Iff \, $\cL$ satisfies interpolation for finitely many sort models.

\noindent
5) In \ref{a32} we can replace clauses (d),(e) by
\mn
\begin{enumerate}
\item[$(d)'$]  for every $\mu < \kappa,\cL$ is $\Delta$-closed in
$\{\lambda:\lambda = \lambda^\mu\}$, see Definition \ref{a31}.
\end{enumerate}
\end{observation}

\begin{PROOF}{\ref{a33}}
Should be clear.
\end{PROOF}

\begin{theorem}
\label{a35}
1) For $\kappa = \beth_\kappa$ the logic 
$\bbL^1_\kappa$ satisfies interpolation.

\noindent
2) Also if $\psi \in \bbL^1_{\le\theta}(\tau_1)$ where $\theta <
\kappa$, the vocabularies $\tau_0 \subseteq \tau_1$ have cardinality 
$\le \theta$ and $\partial = \beth_{\theta^+},\mu = 2^\partial$ 
\then \, we can find a sequence 
$\langle \psi_\alpha:\alpha < (2^\partial)^+\rangle$ of 
members of $\bbL^1_{\le 2^\partial}(\tau_0)$ such that:
\mn
\begin{enumerate}
\item[$(a)$]  $\vdash \psi \rightarrow \psi_\alpha$ for $\alpha < (2^\mu)^+$
\sn
\item[$(b)$]  if $\tau_0 = \tau_1 \cap \tau_2$ and
$\vartheta \in \bbL^1_{\le \theta}(\tau_2)$ a sentence such that $\vdash 
\psi \rightarrow \vartheta$ \then \, $\vdash \psi_\alpha \rightarrow
\vartheta$ for some $\alpha < (2^\mu)^+$.
\end{enumerate}
\end{theorem}

\begin{PROOF}{\ref{a35}}  
1) This means that we should prove the existence of $\psi
\in \bbL^1_\kappa(\tau_0)$ such that $\vdash \psi_1 \rightarrow \psi$
and $\vdash \psi \rightarrow \psi_2$ when we assume:
\mn
\begin{enumerate}
\item[$\boxplus_1$]  $(a) \quad \tau_1 \cap \tau_2 = \tau_0$ 
\sn
\item[${{}}$]  $(b) \quad \psi_\ell \in \bbL^1_\kappa(\tau_\ell) 
\text{ for } \ell=1,2$
\sn
\item[${{}}$]  $(c) \quad \vdash \psi_1 \rightarrow \psi_2$.
\end{enumerate}
\mn
Toward contradiction we assume this fails.
Without loss of generality $\tau_\ell(\ell < 3)$
are vocabularies of cardinality $< \kappa$.

Let $\theta_*$ be such that $\theta_* < \kappa,\theta_* \ge |\tau_1| +
|\tau_2|$ and $\psi_\ell \in \bbL^1_{\le \theta_*}$ for $\ell=1,2$.
For each $\theta \in [\theta_*,\kappa)$ as $\bbL^1_{\le
\theta}(\tau_0)$ is closed under conjunction (and conjunctions) of
$\le \theta$ formulas and negations, clearly
\mn
\begin{enumerate}
\item[$\boxplus_2$]   for every $\alpha < \kappa$ we can choose
$M^1_\alpha,M^2_\alpha$ such that
\mn
\begin{enumerate}
\item[$(a)$]  $M^1_\alpha \text{ is a } \tau_1\text{-model}$
\sn
\item[$(b)$]  $M^2_\alpha \text{ is a } \tau_2\text{-model}$
\sn
\item[$(c)$]  $M^1_\alpha \rest \tau_0,M^2_\theta \rest \tau_0$ are
$\cE^1_{\text{qf}(\tau_0),|\alpha|,\alpha}$-equivalent 
\sn
\item[$(d)$]  $M^1_\alpha \models \psi_1$
\sn
\item[$(e)$]  $M^2_\alpha \models \neg \psi_2$
\end{enumerate}
\end{enumerate}
\mn
and continue naturally as in the proof of \ref{a28}.

\noindent
2) Similarly.
\end{PROOF}

\begin{theorem}
\label{a36}
If $\kappa = \beth_\kappa$ \then \, $\bbL^1_\kappa$ satisfies the addition
and product theorems, that is (all the models are $\tau$-models for a
fixed vocabulary $\tau$; for the sum case we assume $\tau$ has no
function symbols (in particular no individual constants;
or only unary functions) and the relevant models
have disjoint universes):
\mn
\begin{enumerate}
\item[$(a)$]  sum: $\text{\rm Th}_{\bbL^1_\kappa}(M_1 +M_2) = 
\text{\rm Th}_{\bbL^1_\kappa}(M_1) + \text{\rm
Th}_{\bbL^1_\kappa}(M_2)$ that is,
\sn
\item[$(a)'$]  if $\text{\rm Th}_{\bbL^1_\kappa}(M_\ell) = 
\text{\rm Th}_{\bbL^1_\kappa}(N_\ell)$ for $\ell=1,2$ \then \,
$\text{\rm Th}_{\bbL^1_\kappa}(M_1 + M_2) = 
\text{\rm Th}_{\bbL^1_\kappa}(N_1 + N_2)$
\sn
\item[$(b)$]  product - $\text{\rm Th}_{\bbL^1_\kappa}(M_1 \times M_2)
= \text{\rm Th}_{\bbL^1_\kappa}(M_1) \times
\text{\rm Th}_{\bbL^1_\kappa}(M_2)$ that is,
\sn
\item[$(b)'$]  if $\text{\rm Th}_{\bbL^1_\kappa}(M_\ell) = 
\text{\rm Th}_{\bbL^1_\kappa}(N_\ell)$ for $\ell=1,2$ \then \,
$\text{\rm Th}_{\bbL^1_\kappa}(M_1 \times M_2) = \text{\rm
Th}_{\bbL^1_\kappa}(N_1 \times N_2)$
\sn
\item[$(c)$]  moreover, we can replace 
$\text{\rm Th}_{\bbL^1_\kappa}(N)$ by $N
/\cE^1_{\text{\rm qf}(\tau),\theta,\alpha}$ when $\theta <
\kappa,\alpha < \theta^+$.
\end{enumerate}
\end{theorem}

\begin{PROOF}{\ref{a36}}
It suffices to prove clause (c).  We prove it for products.  For the
sums this is easier.
\medskip

\noindent
\underline{Clause (c) for product}:

Clearly it suffices to prove
\mn
\begin{enumerate}
\item[$\boxplus$]  assume $M_\iota = M^\iota_1 \times M^\iota_2$ for
$\iota = 1,2$ are $\tau$-models and $M^1_\ell
\cE^1_{\text{qf}(\tau),\theta,\alpha} M^2_\ell$ for $\ell=1,2$ \then
\, $M_1 \cE^1_{\text{qf}(\tau),\theta,\alpha} M_2$.
\end{enumerate}
\mn
So let $\langle M_{\ell,k}:k \le k(\ell,\iota)\rangle$ be such that
$M_{1,0} = M^1_\ell,M_{2,k} = M^2_\ell$ and $M_{\ell,k}
\cE^0_{\text{qf}(\tau),\theta,\alpha} M_{\ell,k+1}$ for $k <
k(\ell)$.  Let $k(*) = \text{ max}\{k(1),k(2)\}$ and let $M_{\ell,k} =
M_{\ell,k(\ell)}$ if $k(\ell) < k \le k(*)$.  As
$\cE^0_{\text{qf}(\tau),\theta,\alpha}$ is reflexive, still we have
$M_{\ell,k} \cE^0_{\text{qf}(\tau),\alpha} M_{\ell,k+1}$ for
$\ell=1,2$ and $k<k(*)$.

Let $M_k = M_{1,k} \times M_{2,k}$ for $k \le k(*)$ we have to prove
$M_0 \cE^1_{\text{\rm qf}(\tau),\theta,\alpha} M_{k(*)}$; by the
definition of $\cE^1_{\text{\rm qf}(\tau),\theta,\alpha}$ it suffices to
prove for each $k<k(*)$ that
\mn
\begin{enumerate}
\item[$\oplus_k$]  $M_k \cE^0_{\text{qf}(\tau),\theta,\alpha} M_{k+1}$.
\end{enumerate}
\mn
So we have to find a winning strategy for the ISO player for the game
$\Game_{\text{qf}(\tau),\theta,\alpha}[M_k,M_{k+1}]$.

The ISO player restricts itself, in the $n$-th move to 
$(\theta,n)$-states $\bold s_n$ for a play of
$\Game_{\text{qf}(\tau),\theta,\alpha}[M_k,M_{k+1}]$ satisfying,
for some pair $(\bold s_{n,1},\bold s_{n,2})$, that for $\ell=1,2$
\mn
\begin{enumerate}
\item[$(a)$]  $\bold s_{n,\ell}$ is a winning $n$-state for the game
$\Game_{\text{qf}(\tau),\theta,\alpha}[M_{\ell,k},M_{\ell,k+1}]$
\sn
\item[$(b)$]  $A^\ell_{\bold s_n} = A^\ell_{\bold s_{n,1}} \times
A^\ell_{\bold s_{n,2}}$
\sn
\item[$(c)$]   $h_{\bold s_n}((b,c)) = \text{ max}\{h^1_{\bold
s_{n,1}}(b),h^2_{\bold s_{n,\ell}}(c)\}$.
\end{enumerate}
\mn
The rest should be clear.
\end{PROOF}
\bigskip

\centerline {$* \qquad * \qquad *$}

\begin{remark}
\label{a43}
1) Theorem \ref{a36} applies to the monadic version, too.

\noindent
2) Why $\kappa = \beth_\kappa$?  As for any $\theta$ and
$\alpha < \theta^+$ there is a sentence $\psi \in \bbL_{\theta^+}$ and
model $M_\psi$ of $\psi$ of cardinality $\beth_\alpha$ such that any
other model $N$ of $\psi$ can be embedded into $M_\psi$.

\noindent
3) For infinite addition, i.e. $\sum\limits_{s \in I} M^\ell_s$, we have a
   problem: passing to $\cE^0_{\text{qf}(\tau),\theta,\alpha}$ we do
   not know how to uniformize $k(s)$, of course: if $k(s)$ is constant or
   bounded there, then a parallel claim holds, still see \ref{a40} below.

\noindent
4) Even overcoming this obstacle, for infinite product we 
have a problem during a play of
 $\Game_{\text{qf}(\tau),\theta,\alpha}[M,N]$.  The point is about translating
 the $h$'s from here to there and back.  There may well be a
   difference between product (or reduced product) and generalized
product.  Anyhow unlike part (2) this cannot be remedied using
$\cE^0_{\text{qf}(\tau),\theta,\alpha}$ we know that, see \ref{a43}.
\end{remark}

\begin{theorem}
\label{a40}
1) Assume, for simplicity that $\tau$ is a relational vocabulary and
$\langle M^\ell_s:s \in I\rangle$ is a sequence of $\tau$-models with
pairwise disjoint universes, for $\ell=1,2$ and $M_1 = \Sigma\{M^1_s:s
\in S\},M_2 = \Sigma\{M^2_s:s \in S\}$.  If $\kappa = \beth_\kappa$
and $M^1_s,M^2_s$ are $\bbL^1_\kappa$-equivalent
for $s \in S$ \then \, $M_1,M_2$ are $\bbL^1_{\kappa}$-equivalent.

\noindent
2) Moreover for every $\theta,\alpha < \kappa$ there are $\partial,\beta
 < \kappa$ such that for any set $S$ and sequence $\langle M^\ell_s:s
 \in S\rangle$ for $\ell \in \{1,2\}$ of pairwise disjoint
   $\tau$-models we have: if $M^1_s,M^2_s$ are
$\cE_{\text{\rm qf}(\tau),\partial,\beta}$-equivalent for $s \in S$
   then the models $M_1 = \Sigma\{M^1_s:s \in S\},M_2 = 
\Sigma\{M^2_s:s \in S\}$ are $\cE_{\qf(\tau),\theta,\alpha}$-equivalent.

\noindent
3) In part (2) if $M^b_s \cE^0_{\text{\rm qf}(\tau),\theta,\alpha} M^2_j$
for every $s \in S$ \then \, $M_1
\cE^0_{\qf(\tau),\theta,\alpha} M_2$.
\end{theorem}

\begin{PROOF}{\ref{a40}}
1),2) By the proof of \ref{a28}.

\noindent
3) As in the proof of \ref{a36}.
\end{PROOF}

\begin{remark}
\label{a41}
This proof indicates that there are better versions of the game for
which we can choose $(\partial,\beta) = (\theta,\alpha)$.  However, it
is not clear whether a more complicated definition is worth the gain.
\end{remark}
\newpage

\section {$\bbL^1_\kappa$ is strong}

Our first aim is:
\begin{question} 
\label{c2}
How strong is $\bbL^1_\kappa$?  It is more like
$\bbL^{-1}_\kappa$ or $\bbL^0_\kappa$?  E.g. upward LST?  Downward
LST to which cardinals?

So far we have given indication to its being similar to 
$\bbL^{-1}_\kappa$, however we shall show below 
that for the LST, $\bbL^1_\kappa$ is
closer to $\bbL^0_\kappa$ (see \ref{c8}) and fail the theorem on the
theory of an infinite product (which both
$\bbL^0_\kappa,\bbL^{-1}_\kappa$ satisfies).  But restricting
ourselves to cardinals $\lambda = \lambda^{\aleph_0} \ge \kappa$, the
situation for the downward LST is similar to the one for
$\bbL^0_\kappa$, see below.
\end{question}

\begin{definition}
\label{c3}
For a relational vocabulary $\tau$ and cardinal $\theta$ of cofinality
$> \aleph_0$ we say $\bold K$ is a $(\tau,\theta,\omega)$-class of
structures \when \,:
\mn
\begin{enumerate}
\item[$(a)$]  $\bold K$ is a class of $\tau$-models each
of cardinality $\theta$
\sn
\item[$(b)$]  $\bold K$ is closed under submodels of cardinality
$\theta$
\sn
\item[$(c)$]  $\bold K$ is closed under isomorphism
\sn
\item[$(d)$]  $\bold K$ is closed under increasing union of 
$\omega$-chains.
\end{enumerate}
\mn
2) For $\tau,\theta,\bold K$ as above let $\psi_{\bold K} \in
\bbL^0_{(2^\theta)^+,\theta^+}(\tau)$ be a sentence such that for
$\tau$-models $M$ we have: 
$M \models \psi_{\bold K}$ \Iff \, for any $A \subseteq M$ of 
cardinality $\theta$ there is $h:A \rightarrow \omega$ such
that for every large enough $n,M \rest h^{-1}\{n\} \in \bold K$.
\end{definition}

\begin{example}
\label{c3d}
Let $\bold K^{\text{wo}}_\theta = \{(A,<_*):<_*$ is a linear order of $A$ 
and for some sequence $\langle \cU_n:n < \omega\rangle$ we have $(\cU_n,<_*
\rest \, \cU_n) \cong (\theta,<)$ and $\cU_n \subseteq \cU_{n+1}$ and
$A = \cup\{\cU_n:n <\omega\}$.
\end{example}

\begin{claim}
\label{3cp}
For $\tau,\bold K$ as in \ref{c3}.

\noindent
1) The sentence $\psi_{\bold K}$ belongs to
$\bbL_{(2^\theta)^+,\theta^+}(\tau)$ indeed.

\noindent
2) If $\theta + |\bold K / \cong| = \mu$ \then \, moreover 
$\psi_{\bold K} \in \bbL^0_{\mu^+,\theta^+}$.

\noindent
3) Moreover $\psi_{\bold K} \in \bbL^1_{\le \theta}$.
\end{claim}

\begin{PROOF}{\ref{3cp}}
1),2)  Obvious.

\noindent
3) Easy.
\end{PROOF}

\begin{claim}
\label{c4}
1) If $\lambda \ge 2^\theta$ \then \, $\Theta_{\theta,R}$ has a model
of cardinality $\lambda$ iff {\rm cov}$(\lambda,\theta^+,\theta^+,\aleph_0)
= \lambda$, see Remark \ref{c6} below (or \cite[II]{Sh:g}).

\noindent
2) If $\lambda = \lambda^{\aleph_0} \ge 2^{\aleph_0}$ \then \,
   $\Theta_{\theta,R}$ has a model of cardinality $\lambda$ iff
   $\lambda = \lambda^\theta$.
\end{claim}

\begin{PROOF}{\ref{c4}}
1) Read the definitions.

\noindent
2) Check.
\end{PROOF}

\noindent
Complementary is
\begin{claim}
\label{c5}
Assume ($\kappa = \beth_\kappa$ and) $\lambda = 
\text{\rm cov}(\lambda,\kappa,\kappa,\aleph_1)$ and $\tau$ is a
vocabulary of cardinality $< \lambda$.  For every $\tau$-model $M$ of
cardinality $> \lambda$ \underline{there} is $N \subseteq M$ of cardinality
$\lambda$ such that $N \equiv_{\bbL^1_\kappa} M$.

\noindent
2) Assume $\lambda = \text{\rm cov}(\lambda,\beth_\gamma,
\beth_\gamma,\aleph_1)$ for every $\gamma <
\theta^+$.  If $M$ is a $\tau$-model, $|\tau| \le \beth_{\theta^+}$
and $\|M\| > \lambda$ \then \, there is
$N \subseteq M$ of cardinality $\lambda$ such that 
$N \equiv_{\bbL^1_{\le \theta}} M$.
\end{claim}

\begin{PROOF}{\ref{c5}}
Similar to \ref{a21}(4)(c), but in $\omega$ stages deal simultaneously
with $\lambda$ submodels in each stage.
\end{PROOF}

\begin{remark}
\label{c6}
Recall cov$(\lambda,\mu,\theta,\sigma) = \text{ Min}\{|\cP|:\cP
\subseteq [\lambda]^{< \mu}$ and for every $u \in [\lambda]^{<
\theta}$ there are $\alpha < \sigma$ and $u_i \in \cP$ for $i <
\alpha$ such that $u \subseteq \cup\{u_i:i < \alpha\}\}$ for $\lambda
\ge \mu \ge \theta \ge \sigma$.
\end{remark}

\noindent
Concerning the upward LST theorem, the logic $\bbL^1_\kappa$ fail it badly.
\begin{claim}
\label{c8}
There is a sentence $\psi$ from 
$\bbL^1_{\le \aleph_1}$ such that: $\psi$
has a model of cardinality $\lambda$ \Iff \, $\lambda \nrightarrow
(\aleph_1)^{< \omega}_2$.
\end{claim}

\begin{PROOF}{\ref{c8}}
Easy, recalling \ref{3cp}(3).
\end{PROOF}

\begin{discussion}
\label{c10}
We can phrase other relatives.  Also it points to the restrictions when we
are looking for such logics with upward LST.
\end{discussion}

\begin{theorem}
\label{a43}
1) For every $\theta$ there are models $\langle M_n,N_n:n <
\omega\rangle$ with a fixed countable vocabulary
such that $M_n,N_n$ are 
$\bbL^1_{\le \theta}$-equivalent for $n < \omega$ but 
$\prod\limits_{n < \omega} M_n, \prod\limits_{n < \omega} N_n$ 
are not $\bbL^1_{\le 2^{\aleph_0}}$-equivalent.

\noindent
2) Moreover in $\Game_{\qf(\tau),2^{\aleph_0},\omega +1}
 [\prod\limits_{n < \omega} M_n,\prod\limits_{n < \omega} N_n]$ the
   AIS player has a winning strategy.
\end{theorem}

\begin{PROOF}{\ref{a43}}
Let $\tau = \{c,<\} \cup \{F_n:n < \omega\}$ where $F_n$ is a unary function,
$c$ is an individual constant and $<$ is a binary relation.

For every $\alpha$ let $M^n_\alpha$ be the $\tau$-model with universe
$1 + \alpha,c^{M^n_\alpha} = 0,
<^{M^n_\alpha} = \{(i,j):i=0=j$ or $i < j < 1 + \alpha\}$ and 
$F^{M^n_\alpha}_k$ is the identity when $k=n$ and is constantly 0
otherwise.

Fix $\theta \ge 2^{\aleph_0}$, by \ref{a22} for some $\alpha =
\alpha(\theta)$ large enough for each $n$ there is a model $N_n$ which
is $\bbL^1_{\le \theta}$-equivalent to $M^n_\alpha$, but $<^{N_n}
\rest (|N_n| \backslash \{c^{N_n}\})$ is
not well ordered; let $M_n = M^n_\alpha$ and $\langle d_{n,k}:k <
\omega\rangle$ be a $<^{N_n}$-decreasing sequence of members of $N_n$
which are $\ne c^{N_n}$.

Let $M = \prod\limits_{n<\omega} M_n$ and $N = \prod\limits_{n <
\omega} N_n$.  For $\eta \in {}^\omega \omega$ let $d_\eta = \langle
d_{n,\eta(n)}:n < \omega\rangle \in \prod\limits_{n < \omega} N_n =
N$.

Now
\mn
\begin{enumerate}
\item[$\circledast$]  if $\Lambda \subseteq {}^\omega \omega$ is
non-meagre or just not bounded in $({}^\omega
\omega,<_{J^{\text{bd}}_\omega})$ \then \, there is no partial
isomorphism $f$ from $N$ into $M$ with domain $\{d_\eta:\eta \in
\Lambda\}$, i.e. $f$ should preserve all quantifier-free formulas.
\end{enumerate}
\mn
This clearly suffices as if $\Lambda = \cup\{\Lambda_n:n < \omega\}$
then at least for one $n$, the set $\Lambda$ is non-meagre (or just
unbounded in $({}^\omega \omega,<_{J^{\text{bd}}_\omega}))$.

Why is $\circledast$ true?  
Toward contradiction assume that $f$ is such a partial
isomorphism.  By the assumption for some $n(*) < \omega$ the set
$\{\eta(n(*)):\eta \in \Lambda\}$ is infinite, so choose $\bar\eta =
\langle \eta_k:k < \omega\rangle$ such that
\mn
\begin{enumerate}
\item[$(*)_0$]   $\langle\eta_k(n(*)):k <
\omega\rangle$ is strictly increasing in $N_{n(*)}$.
\end{enumerate}
\mn
Now: 
\mn
\begin{enumerate}
\item[$(*)_1$]  if $m(1) < m(2) < \omega$ then for every $n < \omega$ we have
$N_n \models ``F_{n(*)}(d_{n,\eta_{m(2)}(n(*))}) <
F_{n(*)}(d_{n,\eta_{m(1)}(n(*))})"$.
\end{enumerate}
\mn
[Why?  If $n=n(*)$ then $F^{N_n}_{n(*)}$ is the identity and this means
$N_n \models ``d_{n,m(2)} < d_{n,m(1)}"$ which holds by the
choice of $\langle d_{n,k}:k < \omega\rangle$ recalling $(*)_0$.  
If $n < \omega \wedge n \ne n(*)$ then $F^{N_n}_{n(*)}$ is 
constantly 0 so this means $N_n
\models ``0 < 0"$ which holds by the choice of $N_n$ (and
$<^{M^n_\alpha}$), so we are done.]
\mn
\begin{enumerate}
\item[$(*)_2$]  $N \models ``F_{n(*)}(d_{\eta_{m(2)}}) <
F_{n(*)}(d_{\eta_{m(1)}})"$ for $m(1) < m(2) < \omega$.
\end{enumerate}
\mn
[Why?  By $(*)_1$ and the definition of product.]

Also
\mn
\begin{enumerate}
\item[$(*)_3$]  if $m < \omega$ then
$N_{n(*)} \models ``c \ne F_{n(*)}(d_{n(*),\eta_m(n(*))})"$
hence $N \models ``c \ne F_{n(*)}(d_{\eta_n})"$
\sn
\item[$(*)_4$]  let $\nu_k = f(d_{\eta_k})$ so $\nu_k = \langle a_{k,n}:n
< \omega\rangle \in \prod\limits_{n < \omega} M_n$.
\end{enumerate}
\mn
But $f$ is a partial isomorphism so by $(*)_2 + (*)_3 + (*)_4$ we have
\mn
\begin{enumerate}
\item[$(*)_5$]  $M \models ``F_{n(*)}(\nu_{m(2)}) <
F_{n(*)}(\nu_{m(1)}) \ne c"$ for $m(1) < m(2) < \omega$.
\end{enumerate}
\mn
Hence for every $n,M_n \models ``F_{n(*)}(\nu_{m(2)}(n)) <
F_{n(*)}(\nu_{m(1)}(n))"$ for $m(1) < m(2) < \omega$.  Also by the
choice of $F^{M_n}_{n(*)}$ we have $n \ne n(*)
\Rightarrow M_n \models ``c = F_{n(*)}(\nu_{m(1)}(n))"$ but $M \models
``c \ne F_{n(*)}(\nu_{m(1)})"$ hence $M_{n(*)}
\models ``c \ne F_{n(*)}(\nu_{m(1)}(n(*)))"$, i.e. $\nu_{m(1)}(n(*)) \ne
0$.  Together $\langle \nu_m(n(*)):m < \omega\rangle$ is
$<^{M_n}$-decreasing in $M_{n(*)} \backslash \{0\}$ contradiction.
\end{PROOF}


\end{document}